\documentclass[a4paper,12pt]{article}
\usepackage{amsmath}
\usepackage[mathscr]{eucal}
\usepackage{amssymb}
\usepackage{latexsym}
\usepackage{amsthm}
\usepackage{amscd}
\theoremstyle{plain}
\newtheorem{theorem}{Theorem}

\newtheorem{lemma}{Lemma}[section]
\newtheorem{proposition}{Proposition}

\theoremstyle{remark}
\newtheorem{remark}{Remark}
\theoremstyle{definition}

\usepackage[dvips]{graphicx}
\usepackage{psfrag}
\begin{document}
\title{ Notes on hyperelliptic fibrations of genus $3$, II \\ 
\footnotetext{During this study, the author was supported 
by DFG Forshergruppe 790
``Classification of algebraic surfaces and compact complex manifolds''.}
}
\author{Masaaki Murakami \footnotetext{During the execution of 
this study, the author was supported by...}}
\date{}
\maketitle
\begin{abstract}
This is the part II of the series under the same title. 
In part I, 
using the approach developed by 
Catanese--Pignatelli \cite{fibrationsI'}, we gave a structure 
theorem for hyperelliptic genus $3$ fibrations 
all of whose fibers are $2$--connected (\cite{notesonI}).  
In this part II, we shall give a structure theorem for 
smooth deformation families of these  
to non-hyperelliptic genus $3$ fibrations.  
As an application, we shall give a set of sufficient conditions 
for our genus $3$ hyperelliptic fibration above to allow deformation  
to non-hyperelliptic fibrations, and use this to study certain 
minimal regular surfaces with $c_1^2 = 8$ and $p_g =4$:
we shall find in the moduli space 
two strata $\mathcal{M}_0^{\sharp} $ and $\mathcal{M}_0^{\flat}$ 
(each of dimensions $32$ and $30$, respectively), 
and show that Bauer--Pingatelli's stratum $\mathcal{M}_0$ (\cite{caninvpg4c8}) 
and its $26$-dimensional substratum are at the boundary 
of these new strata $\mathcal{M}_0^{\sharp} $ and $\mathcal{M}_0^{\flat}$,
respectively.   
\end{abstract}

\section{Introduction}  \label{sectn:introduction}

In classification of algebraic surfaces of general type, 
we sometimes encounter a situation where we have classified all the 
possible types of the surfaces of a class, but 
not yet having clarified the global structure of moduli space.  
Such situation appears for example in the case of 
minimal regular surfaces with $c_1^2 = 6$ and $p_g = 4$, 
where Horikawa classified the surfaces 
into $11$ types, and showed, gluing these strata, that 
their moduli space has at most $3$ connected components (\cite{smallc1-3}). 
Although Bauer--Catanese--Pignatelli \cite{bcpmodulic6} later 
succeeded in showing that the moduli space has at most two 
connected components, the determination of 
the number of the connected components still remains 
as a challenging problem. 
Since we have such cases, 
it is of interest to develop a tool with which we can glue 
strata of a moduli space.

In the present paper, as an example of such attempts, 
we shall give a structure theorem for smooth 
non-hyperelliptic deformation families 
(over smooth base spaces) of 
the genus $3$ hyperelliptic fibrations treated in \cite{notesonI}. 
Here by a non-hyperelliptic deformation family we mean a 
deformation family 
(of fibrations) 
whose general fibers give 
non-hyperelliptic fibrations.  
We shall also give using this structure theorem  
a set of sufficient conditions 
for our hyperelliptic genus $3$ fibration to allow a 
non--hyperelliptic deformation (Theorem \ref{thm:exstcthm}),  
and then 
use it to study certain classes of minimal regular surfaces 
with $p_g = 4$ and $c_1^2  = 8$.   
Note that, for our purpose,  
the explicit construction in the proof of Theorem \ref{thm:exstcthm} 
is as important as 
the assertion of the theorem itself.  

Our structure theorem is, just like those given in \cite{fibrationsI'}, 
a one--to--one correspondence between the deformation families 
and admissible $5$-tuples. 
Since outline of the proof is similar to that for 
case of genus $3$ fibrations in \cite{fibrationsI'}, the second part
(i.e., the set of sufficient conditions and its application to 
the study of minimal regular surfaces with $p_g = 4$ and $c_1^2 = 8$)
in a sense is the main part of the present paper. 

The study of minimal regular surfaces with $p_g = 4$ has  
long history from Enriquess \cite{enriquessuperfici},  
and since then 
there have been a lot of works on the surfaces with 
$p_g = 4$. After the complete classification 
of the case $c_1^2 = 7$ by I.\,Bauer \cite{pg4c7bauer}, 
the next case $c_1^2 = 8$ has been of interests of 
albebro-geometers working on algebraic surfaces
(see, e.g., \cite{supinopg4c8}, \cite{oliverioevenpg4c8}, 
\cite{caninvpg4c8}, \cite{catlipign}). 
In the present paper we shall find two strata  
$\mathcal{M}_0^{\sharp}$ and $\mathcal{M}_0^{\flat}$ 
(in the moduli space) of respective dimensions  
$\dim \mathcal{M}_0^{\sharp} = 32$ and 
$\dim \mathcal{M}_0^{\flat} = 30$.  
Then as an application of the sufficient conditions, 
we shall show that 
the stratum 
$\mathcal{M}_0$ in \cite{caninvpg4c8} and its $26$-dimensional substratum  
are at the boundary of our two strata  
$\mathcal{M}_0^{\sharp}$ and $\mathcal{M}_0^{\flat}$, 
respectively.
We remark here that, for both of our new strata  
$\mathcal{M}_0^{\sharp}$ and $\mathcal{M}_0^{\flat}$,  
general members have no canonical involution. 
It is likely that members of 
$\mathcal{M}_0^{\sharp}$ belong to a family 
constructed by Ciliberto in \cite{cilcansfpg4}, but as 
boundary cases (see Remark \ref{rm:refmsharp}). 
As for members of $\mathcal{M}_0^{\flat}$, 
to the authors knowledge, they are new (see Remark \ref{rm:refmflat}). 

The present paper is organized as follows.    
In section \ref{scn:strthm}, we shall prove the existence of a  
natural one--to--one correspondence between the isomorphism classes of  
non-hyperelliptic deformation families (of a $2$--connected 
hyperelliptic genus $3$ fibration) and the isomorphism classes of  
admissible $5$-tuples. 
Although outline of the proof is similar to that in 
\cite{fibrationsI'}, we replace details by more concrete computations. 
This is because at the last part of the proof, because of 
the lack of compactness of the total space of the deformation 
families, we want to execute computations using the differential 
forms in stead of cohomological definition of the dualizing sheaves.
These concrete descriptions are implicitly used also in the 
proof of Theorem \ref{thm:exstcthm}.     
 
In section \ref{scn:applsfswith}, 
we shall give a set of sufficient conditions 
(Theorem \ref{thm:exstcthm}), and apply it to the study of 
minimal regular surfaces with $p_g = 4$ and $c_1^2 = 8$. 
The proof of Theorem \ref{thm:exstcthm} 
is based on the results in Section \ref{scn:strthm}. 
As for the application, we shall first find the strata  
$\mathcal{M}_0^{\sharp}$ and $\mathcal{M}_0^{\flat}$, then 
study general members of Bauer-Pignatelli's $\mathcal{M}_0$ and 
its substratum, and finally deform the members of 
$\mathcal{M}_0^{\sharp}$ 
to those of $\mathcal{M}_0^{\sharp}$ and $\mathcal{M}_0^{\flat}$. 
\medskip

{\sc Acknowledgment}   
 \medskip

\noindent
The author expresses deepest gratitude to Prof.\ Fabrizio Catanese 
and Prof.\ Ingrid Bauer 
for the comfortable environment, financial support, 
and discussions. 
During this study, the author was supported by DFG Forshergruppe 790
``Classification of algebraic surfaces and compact complex manifolds''.  

\medskip

{\sc Notation}
\medskip

\noindent
In this article, the symbol $k$ always denotes the complex 
number field $\mathbb{C}$.  
If $V$ is a locally free sheaf of finite rank on a scheme, 
$\mathrm{rk} \, V$ denotes its rank. 
The symbols $\mathbb{S} (V)$ and $\mathbb{S}^n (V)$ denote 
the symmetric tensor algebra (of $V$) and its homogeneous part 
of degree $n$, respectively. 
If $R$ is a graded algebra, $R_j$ denotes its homogeneous part 
of degree $j$. Thus, for example, for the polynomial ring 
$k[x_0, x_1, x_2]$ over the complex number field $k= \mathbb{C}$, 
an element in $k[x_0, x_1, x_2]_j$ is 
a homogeneous polynomial in $x_0$, $x_1$, $x_2$ of degree $j$. 
The symbol $\amalg$ means taking the disjoint union of sets. 
If $p$ is a point of a scheme, $k(p)$ denotes the residue field 
at $p$ of this scheme. If $\mathcal{F}$ is a sheaf on a scheme, 
$\mathcal{F}_p$ denotes the stalk at $p$ of $\mathcal{F}$, and  
$\mathrm{supp}\, \mathcal{F}$, the support of $\mathcal{F}$. 
If moreover the scheme is over $k$, and $\mathcal{F}$ is 
coherent, $h^i (\mathcal{F})$ denotes the dimension over $k$ 
of the $i$-th cohomology group $H^i (\mathcal{F})$ of $\mathcal{F}$.
If $\sigma : \mathcal{F} \to \mathcal{G}$ is a morphism of sheaves, 
$\mathrm{Im}\, \sigma$ denotes the image of $\sigma$. 
The identity morphism of $\mathcal{F}$ is denoted by 
$\mathrm{id}_{\mathcal{F}}$. 
If $W$ is a subscheme, $\mathcal{F} |_W$ denotes the restriction 
of $\mathcal{F}$ to $W$.   
If $S$ is a scheme, and $D$, a Cartier divisor on $S$, then 
$\mathcal{O}_S$ and $\mathcal{O}_S (D)$ denote the structure sheaf 
of $S$ and the invertible sheaf associated to $D$, respectively. 
If $S$ is a smooth variety over $k$, the symbol $K_S$ 
as usual denotes the canonical divisor of $S$. 
If $\mathcal{M}$ is a moduli space, and $S$, a member of $\mathcal{M}$, 
then $[S] \in \mathcal{M}$ denotes the point of $\mathcal{M}$ corresponding 
to the isomorphism class of $S$.   

\section{Structure theorem}   \label{scn:strthm}

Let us first fix the terminology. 
Assume we are given a non-singular variety $T$, 
a proper smooth family $\pi_{\mathcal{B}} : \mathcal{B} \to T$ 
of irreducible curves of genus $\frak{b}$,  
a proper sooth family $\pi_{\mathcal{S}}: \mathcal{S} \to T$ 
of irreducible surfaces, and 
a proper flat surjective morphism $\frak{f}: \mathcal{S} \to \mathcal{B}$
such that $\pi_{\mathcal{S}} = \pi_{\mathcal{B}} \circ \frak{f}$.         
If we fix a closed point $t_0 \in T$, 
we shall call the collection of data 
$(\mathcal{S}, \mathcal{B}, T, t_0, 
\frak{f}, \pi_{\mathcal{S}},  \pi_{\mathcal{B}})$, 
a family of deformations of the fibration 
$f: S = \pi_{\mathcal{S}}^{-1} (t_0) \to B= \pi_{\mathcal{B}}^{-1} (t_0)$ 
(where $f$ is the restriction of $\frak{f}$). 
For each closed point $t \in T$, we denote by 
$f_t : S_t = \pi_{\mathcal{S}}^{-1} (t) \to B_t= \pi_{\mathcal{B}}^{-1} (t)$ 
the restriction of $\frak{f}$ to the fibers $S_t$ and $B_t$ over $t$.

In the present paper, we study the case where the central 
fibration $f: S \to B$ is a genus $3$ relatively minimal 
hyperelliptic fibration all of whose fibers are $2$-connected, and such that 
for general  $t \in T$ the fibration $f_t : S_t \to B_t$ is 
non-hyperelliptic. In this case, if we replace the parameter space 
$T$ by smaller one, we may assume that the $2$-connectedness property 
holds for any $t \in T$. So in the present paper, we assume that 
for any $t \in T$, the $f_t : S_t \to B_t$ is a relatively minimal  
genus $3$ fibration all of whose fibers are $2$-connected.   
We shall call such deformation family a ($2$--connected) non-hyperelliptic 
deformation family of the $2$-connected hyperelliptic genus $3$ 
fibration $f: S \to B$. 

Given such a family, we define the coherent sheaf 
$\omega_{\mathcal{S}| \mathcal{B}}$ by 
$\omega_{\mathcal{S}| \mathcal{B}} 
= \mathcal{O}_{\mathcal{S}} (K_{\mathcal{S}} - \frak{f}^* K_{\mathcal{B}})$, and 
denote by 
$\mathcal{V}_n = \mathfrak{f}_* (\omega_{\mathcal{S}| \mathcal{B}}^{\otimes n}) $ 
the direct image by $\mathfrak{f}$ of its $n$-th tensor power 
$\omega_{\mathcal{S}| \mathcal{B}}^{\otimes n}$. 
Then for any natural number $n \geq 0$, the sheaf $\mathcal{V}_n$ is 
locally free, and its rank is given by  
$\mathrm{rk}\, \mathcal{V}_n = 4n -2$ for $n \geq 2$, and 
$\mathrm{rk}\, \mathcal{V}_n = 3$ for $n =1$. 
We define the relative canonical algebra 
$\mathcal{R}$ by $\mathcal{R} = \mathcal{R} (\frak{f})
= \oplus_{n=0}^{\infty} \mathcal{V}_n$, and 
the relative canonical model $\mathcal{X}$ 
by $\frak{g} : \mathcal{X} = \mathcal{P}\mathrm{roj}\, \mathcal{R} 
\to \mathcal{B}$, where $\frak{g}$ is the structure morphism. 
We put $\pi_{\mathcal{X}} = \pi_{\mathcal{B}} \circ \frak{g} 
: \mathcal{X} \to \mathcal{B}$. 
Then it is easy to see that for any closed point $t \in T$ 
the natural restriction $g_t : X_{t} = \pi_{\mathcal{X}}^{-1} (t) \to 
B_t = \pi_{\mathcal{B}}^{-1} (t)$ is the relative 
canonical model of $f_t : S_{t} \to B_{t}$. 

As in Catanese--Pignatelli \cite{fibrationsI'}, we consider the 
natural morphism 
$\bar{\sigma}_n : \mathbb{S}^n (\mathcal{V}_1) \to \mathcal{V}_n$ 
determined by the multiplication structure of the relative canonical 
algebra $\mathcal{R}$. 
Let us denote by $\bar{\mathcal{T}}_n = \mathrm{Cok}\, \bar{\sigma}_n$
the cokernel of the morphism $\bar{\sigma}_n$. 

By the same method as in Catanese--Pignatelli \cite{fibrationsI'}, 
we can prove the following:

\begin{lemma}  \label{lm:localstrrelcan}
For any closed point $b \in \mathcal{B}$, 
there exists an affine neighbourhood 
$\mathcal{U} \subset \mathcal{B}$ of $b$ 
such that $\mathcal{R} (\mathcal{U})$ is of the form 
\[
\mathcal{R} (\mathcal{U}) \simeq 
\mathcal{O}_{\mathcal{B}} (\mathcal{U}) [x_0, x_1, x_2, y]
/ (hy - \bar{Q} (x_0, x_1, x_2), \, y^2 - \bar{P} (x_0, x_1, x_2))
\]
$($$\deg x_i = 1$ for $0 \leq i \leq 2$ and $\deg y = 2$$)$ 
as graded $\mathcal{O}_{\mathcal{B}} (\mathcal{U})$--module, 
where 
$\bar{Q} (x_0, x_1, x_2) \in 
\mathcal{O}_{\mathcal{B}} (\mathcal{U}) [x_0, x_1, x_2]_2$ and 
$\bar{P} (x_0, x_1, x_2) \in 
\mathcal{O}_{\mathcal{B}} (\mathcal{U}) [x_0, x_1, x_2]_4$ are 
homogeneous polynomials in $x_0$, $x_1$, $x_2$ of degree 
$2$ and $4$, respectively,  
and 
$h \neq 0 \in \mathcal{O}_{\mathcal{B}} (\mathcal{U})$ 
is a non--zero function on $\mathcal{U}$ that identically vanishes 
on $B_{t_0} \cap \mathcal{U}$. 
\end{lemma}

We remark here that if we do not take the completion of the 
square with respect to $y$, we can set two relations in 
the relative canonical algebra also in the following form:  
$hy - x_2^2 + \bar{Q}_1 x_2 + \bar{Q}_2$ and 
$y^2 - \bar{P_1} x_2 y - \bar{P_2} y - \bar{P_3} x_2 - \bar{P_4}$, 
where
$\bar{P_i} \in {O}_{\mathcal{B}} (\mathcal{U}) [x_0, x_1]_i$ 
and $\bar{Q}_j \in {O}_{\mathcal{B}} (\mathcal{U}) [x_0, x_1]_j$, 
and such that the restriction $\bar{P_1} |_{B_t \cap \mathcal{U}}$ 
and $\bar{P_2} |_{B_t \cap \mathcal{U}}$ to $B_t \cap \mathcal{U}$ 
identically vanish.

The functions $h$'s in the lemma above 
patch to define an effective divisor $\bar{\tau}$
on $\mathcal{B}$,  
and with this $\bar{\tau}$ 
the sheaf $\bar{\mathcal{T}} = \bar{\mathcal{T}}_2$ is
an invertible $\mathcal{O}_{\bar{\tau}}$--module. 
Note that the fibration $f_t : S_t \to B_t$ is hyperelliptic 
if and only if $B_t \subset \mathrm{supp} \bar{\tau}$.

\begin{lemma}   \label{lm:rkimsigman}
For any integer $n \geq 2$, the image $\mathrm{Im}\, \bar{\sigma}_n$
is a locally free sheaf of rank $4n -2$. 
The sheaf $\bar{\mathcal{T}}_n$ is a locally free 
$\mathcal{O}_{\bar{\tau}}$--module of rank $2n - 3$. 
\end{lemma} 

Proof. 
We only need to check the assertion in a neighbourhood of each 
$b \in \mathrm{supp}\, \bar{\tau}$. 
Take the sections $x_i$'s and $y$ such that the relations in the 
relative canonical algebra are of the form
$hy - x_2^2 + \bar{Q}_1 x_2 + \bar{Q}_2$ and 
$y^2 - \bar{P_1} x_2 y - \bar{P_2} y - \bar{P_3} x_2 - \bar{P_4}$
and such that $\bar{P_1} |_{B_t \cap \mathcal{U}} \equiv
\bar{P_2} |_{B_t \cap \mathcal{U}} \equiv 0$. 
For each $n \geq 0$ we denote by $(*)_n$ the set of monomials 
in $x_0$, $x_1$ of degree $n$.   
Then the disjoint union 
$(*)_n \amalg  ((*)_{n-1} \cdot x_2 ) 
\amalg ((*)_{n-2}\cdot  y) \amalg ((*)_{n-3} \cdot x_2 y)$ forms a 
base of the free $\mathcal{O}_{\mathcal{B}} (\mathcal{U})$--module  
$\mathcal{V}_n (\mathcal{U})$, and by the two relations above 
we see easily that the module 
$(\mathrm{Im}\, \bar{\sigma}_n) (\mathcal{U})$ is generated in 
$\mathcal{V}_n (\mathcal{U})$ by the disjoint union 
$(*)_n \amalg ((*)_{n-1}\cdot x_2) 
\amalg ((*)_{n-2}\cdot hy) \amalg ((*)_{n-3}\cdot hx_2 y)$.   
Then the assertion easily follows. \qed

\begin{lemma}  \label{lm:lnandsgmn}
Let $\mathcal{L}_n = \ker \bar{\sigma}_n$ be the kernel of the morphism 
$\bar{\sigma}_n : \mathbb{S}^n (\mathcal{V}_1) \to \mathcal{V}_n$. 
Then the following hold:

$1$$)$ $\mathcal{L}_2 \simeq \mathcal{L}_3 \simeq 0$. 
Moreover for each $n \geq 4$, the sheaf $\mathcal{L}_n$ is 
locally free of rank $(n-2)(n-3)/2$. In particular 
$\mathcal{L}_4$ is an invertible sheaf. 

$2$$)$ For each $n \geq 4$, the natural morphism 
$\mathcal{L}_4 \otimes \mathbb{S}^{n-4} (\mathcal{V}_1) 
\to \mathcal{L}_n$ is an isomorphism. 
In particular the relative $1$--canonical image 
$\varSigma$ is a divisor in $\mathbb{P} (\mathcal{V}_1)$ 
belonging to the linear system 
$| \mathcal{O}_{\mathbb{P} (\mathcal{V}_1)} (4) \otimes 
\mathrm{pr}_{\mathbb{P}(\mathcal{V}_1)}^* \mathcal{L}_4^{\otimes (-1)} |$, 
where $\mathrm{pr}_{\mathbb{P}(\mathcal{V}_1)} :
 \mathbb{P}(\mathcal{V}_1) \to \mathcal{B}$ is the structure morphism.    
\end{lemma}

Proof. 
The assertion $1$) follows from the exact sequence 
$0 \to \mathcal{L}_n \to \mathbb{S}^n (\mathcal{V}_1) 
\to \mathrm{Im}\, \bar{\sigma}_n \to 0$, since 
$\mathbb{S}^n (\mathcal{V}_1)$ and 
$\mathrm{Im} \bar{\sigma}_n $ are locally free sheaves 
of rank
$\mathrm{rk}\, \mathbb{S}^n (\mathcal{V}_1) = (n+2)(n+1)/2$ and 
$\mathrm{rk}\, \mathrm{Im} \bar{\sigma}_n = 4n -2$, respectively. 
Let us prove the assertion $2$). 
Since $\mathrm{Im}\, \bar{\sigma}_4$ is locally free, we see that 
for any closed point $b \in \mathcal{B}$ the natural morphism 
$\mathcal{L}_4 \otimes k(b) \to \mathbb{S}^4 (\mathcal{V}_1) \otimes k(b)$
is injective. Thus 
for any $n \geq 4$ the natural morphism 
$(\mathcal{L}_4 \otimes \mathbb{S}^{n-4} (\mathcal{V}_1)) \otimes k(b) 
\to \mathbb{S}^n (\mathcal{V}_1) \otimes k(b) $ is injective, which implies 
also the injectivity of the morphism
$(\mathcal{L}_4 \otimes \mathbb{S}^{n-4} (\mathcal{V}_1)) \otimes k(b) 
\to \mathcal{L}_n \otimes k(b)$. 
This implies the assertion, since we have 
$\mathrm{rk}\, (\mathcal{L}_4 \otimes \mathbb{S}^{n-4} (\mathcal{V}_1))
= \mathrm{rk}\, \mathcal{L}_n = (n-2)(n-3)/2$. \qed

As in Catanese--Pignatelli \cite{fibrationsI'}, we consider the 
morphism $\bar{\frak{c}} : \mathbb{S}^2 ( \bigwedge^2 \mathcal{V}_1) \to 
\mathbb{S}^2 (\mathbb{S}^2 (\mathcal{V}_1))$ 
given by 
$(a \wedge b)(c \wedge d) \mapsto (ac)(bd) - (ad)(bc)$, 
and define the coherent sheaf $\tilde{\mathcal{V}}_4$ 
on $\mathcal{B}$ by 
$\tilde{\mathcal{V}}_4 
= \mathrm{Cok}\, (\mathbb{S}^2 (\bar{\sigma}_2) \circ \bar{\frak{c}} 
: \mathbb{S}^2 ( \bigwedge^2 \mathcal{V}_1) \to 
\mathbb{S}^2 (\mathcal{V}_2))$.  
By the same argument as in  
\cite[Lemma 7.6]{fibrationsI'}, we see that 
$\tilde{\mathcal{V}}_4$ is a locally free sheaf of rank 
$\mathrm{rk}\, \tilde{\mathcal{V}}_4 = 15$.
Note that we have a natural short exact sequence 
\[
0 \to \mathbb{S}^2 ( \bigwedge^2 \mathcal{V}_1) 
  \to \mathbb{S}^2 (\mathbb{S}^2 (\mathcal{V}_1))
  \to \mathbb{S}^4 (\mathcal{V}_1) \simeq \mathrm{Cok}\, \bar{\frak{c}}
  \to 0.
\]

\begin{lemma}  \label{lm:veroimgradealg}
Let $b \in \mathcal{B}$ be a closed point. 
Take an affine neighbourhood $\mathcal{U} \subset \mathcal{B}$ of $b$ 
as in Lemma \ref{lm:localstrrelcan}, i.e., such that 
$\mathcal{R} (\mathcal{U}) \simeq 
\mathcal{O}_{\mathcal{B}} (\mathcal{U}) [x_0, x_1, x_2, y]/
(hy - \bar{Q}, y^2 - \bar{P})$ for a suitable choice of 
$x_i$'s  $ \in \mathcal{V}_1 (\mathcal{U})$ and 
$y \in \mathcal{V}_2 (\mathcal{U})$, and of   
$\bar{Q} \in \mathcal{O}_{\mathcal{B}} (\mathcal{U}) [x_0, x_1, x_2]_2$
and 
$\bar{P} \in \mathcal{O}_{\mathcal{B}} (\mathcal{U}) [x_0, x_1, x_2]_4$.
Define the graded algebra $\tilde{\mathcal{R}} (\mathcal{U})$ by  
$\tilde{\mathcal{R}} (\mathcal{U}) = 
\mathcal{O}_{\mathcal{B}} (\mathcal{U}) [x_0, x_1, x_2, y]
/(hy - \bar{Q})$. 
Then there exist natural isomorphisms 
$\tilde{\mathcal{R}} (\mathcal{U})_2 \simeq \mathcal{V}_2 (\mathcal{U})$
and 
$\tilde{\mathcal{R}} (\mathcal{U})_4 \simeq \tilde{\mathcal{V}}_4 (\mathcal{U})$
for which the morphism $\mathbb{S}^2 (\mathcal{V}_2) (\mathcal{U}) \to 
\tilde{\mathcal{V}}_4 (\mathcal{U})$ induced by the 
multiplication morphism $\mathbb{S}^2 (\tilde{\mathcal{R}} (\mathcal{U})_2) 
\to \tilde{\mathcal{R}} (\mathcal{U})_4$ coincides 
with the natural projection $\mathbb{S}^2 (\mathcal{V}_2) (\mathcal{U}) \to 
\tilde{\mathcal{V}}_4 (\mathcal{U})$ determined by the definition of 
$\tilde{\mathcal{V}}_4
= \mathrm{Cok}\, (\mathbb{S}^2 (\bar{\sigma}_2) \circ \bar{\frak{c}} 
: \mathbb{S}^2 ( \bigwedge^2 \mathcal{V}_1) \to 
\mathbb{S}^2 (\mathcal{V}_2))$.
\end{lemma}

Proof. 
Let  
$\tilde{\mathcal{R}} (\mathcal{U}) \to \mathcal{R} (\mathcal{U})$ 
be the natural projection. 
Then via this projection we have isomorphisms 
$\tilde{\mathcal{R}} (\mathcal{U})_i \simeq \mathcal{R} (\mathcal{U})_i$
for $i = 1$, $2$. 
Consider the multiplication morphism
\begin{equation}   \label{eqn:multmor}
\mathbb{S}^2(\tilde{\mathcal{R}} (\mathcal{U})_2) \to 
\tilde{\mathcal{R}} (\mathcal{U})_4,
\end{equation}
and denote by 
$\mathbb{S}^2 (\mathbb{S}^2(\tilde{\mathcal{R}} (\mathcal{U})_1))
\to \tilde{\mathcal{R}} (\mathcal{U})_4$ the 
composite of two morphisms 
$\mathbb{S}^2 (\mathbb{S}^2(\tilde{\mathcal{R}} (\mathcal{U})_1))
\to \mathbb{S}^2(\tilde{\mathcal{R}} (\mathcal{U})_2)$ and 
(\ref{eqn:multmor}), where 
the $\mathbb{S}^2 (\mathbb{S}^2(\tilde{\mathcal{R}} (\mathcal{U})_1))
\to \mathbb{S}^2(\tilde{\mathcal{R}} (\mathcal{U})_2)$ is the 
morphism induced by the multiplication 
$\mathbb{S}^2(\tilde{\mathcal{R}} (\mathcal{U})_1)
\to \tilde{\mathcal{R}} (\mathcal{U})_2$. 
Consider also the the composite 
$\mathbb{S}^2 (\bigwedge^2 \tilde{R} (\mathcal{U})_1)
\to \tilde{R} (\mathcal{U})_4$
of the two morphisms
\[
\mathbb{S}^2 (\bar{\sigma}_2) \circ \bar{\frak{c}} : 
\mathbb{S}^2 (\bigwedge^2 \tilde{\mathcal{R}} (\mathcal{U})_1)
(\simeq
\mathbb{S}^2 (\bigwedge^2 \mathcal{V}_1 (\mathcal{U})))
\to
\mathbb{S}^2 (\tilde{\mathcal{R}} (\mathcal{U})_2)
(\simeq
\mathbb{S}^2 (\mathcal{V}_2 (\mathcal{U})))
\]
and (\ref{eqn:multmor}). 
Since the 
$\mathbb{S}^2 (\mathbb{S}^2(\tilde{\mathcal{R}} (\mathcal{U})_1))
\to \tilde{\mathcal{R}} (\mathcal{U})_4$ above factors through 
the natural multiplication morphism 
$\mathbb{S}^4(\tilde{\mathcal{R}} (\mathcal{U})_1)
\to \tilde{\mathcal{R}} (\mathcal{U})_4$, 
we see that the composite 
$\mathbb{S}^2 (\bigwedge^2 \tilde{R} (\mathcal{U})_1)
\to \tilde{R} (\mathcal{U})_4$ above is a zero morphism. 
This implies that the morphism (\ref{eqn:multmor})
induces a surjection 
$\tilde{\mathcal{V}}_4 (\mathcal{U}) 
\simeq \mathrm{Cok}\, (
\mathbb{S}^2 (\bar{\sigma}_2) \circ \bar{\frak{c}} )
\to
\tilde{R}(\mathcal{U})_4$.
Since both 
$\tilde{\mathcal{V}}_4 (\mathcal{U}) $ and  
$ \tilde{R}(\mathcal{U})_4$ 
are free of rank $15$, we have 
the assertion. \qed

Since the natural morphism 
$\mathbb{S}^2(\mathbb{S}^2 (\mathcal{V}_1)) 
(\to \mathbb{S}^2 (\mathcal{V}_2))
\to \mathcal{V}_4$ factors through 
$\mathbb{S}^2(\mathbb{S}^2 (\mathcal{V}_1)) \to \mathbb{S}^4 (\mathcal{V}_1)
\simeq \mathrm{Cok}\, \bar{\frak{c}}$, 
this $\mathbb{S}^2 (\mathcal{V}_2) \to \mathcal{V}_4$ induces 
the natural projection 
$\tilde{\mathcal{V}}_4 \to \mathcal{V}_4$.
We denote by $\mathcal{L}_4^{\prime}$ the kernel of this 
induced projection $\tilde{\mathcal{V}}_4 \to \mathcal{V}_4$.
Note that we have a natural morphism 
$\mathbb{S}^4 (\mathcal{V}_1) \simeq 
\mathrm{Cok}\, \bar{\frak{c}} \to 
\tilde{\mathcal{V}}_4 = 
\mathrm{Cok}\, (\mathbb{S}^2 (\bar{\sigma_2}) \circ \bar{\frak{c}})$. 
Since the composite  
$\mathbb{S}^4 (\mathcal{V}_1) \to \mathcal{V}_4$ of this 
morphism with the projection 
$\tilde{\mathcal{V}}_4 \to \mathcal{V}_4$ 
coincides with the multiplication 
$\bar{\sigma_4} : \mathbb{S}^4 (\mathcal{V}_1) \to \tilde{\mathcal{V}}_4$,   
we obtain the 
commutative diagram
\begin{equation*} 
\begin{CD}  
 0  @>\text{$$}>> \mathcal{L}_4  @>\text{$$}>> \mathbb{S}^4(\mathcal{V}_1)  @>\text{$$}>> 
   \mathbb{S}^4(\mathcal{V}_1) / \mathcal{L}_4 @>\text{$$}>> 0\\
  @VV\text{$$}V  @VV\text{$$}V  @VV\text{$$}V  @VV\text{$$}V @VV\text{$$}V\\
  0  @>\text{$$}>> \mathcal{L}_4^{\prime} @>\text{$$}>> \tilde{\mathcal{V}}_4 @>\text{$$}>> \mathcal{V}_4  @>\text{$$}>> 0 \, , 
\end{CD}
\end{equation*}
where the two rows are exact. 

\begin{lemma}   \label{lm:xhpsfofw}
Let $\mathcal{W} \subset \mathbb{P} (\mathcal{V}_2)$ be the image 
of the rational veronese map 
$\mathbb{P} (\mathcal{V}_1) - - \to \mathbb{P} (\mathcal{V}_2)$. 
Then the relative canonical model $\mathcal{X}$ of 
$\mathcal{S}$ is a hypersurface of $\mathcal{W}$ 
cut out by $\mathcal{L}_4^{\prime}$. 
\end{lemma}

Proof. 
Let $b \in \mathrm{supp}\, \bar{\tau}$, and take 
an affine neighbourhood $\mathcal{U} \subset \mathcal{B}$ of $b$
as in Lemma \ref{lm:localstrrelcan}.  
Consider the graded algebra 
$\tilde{\mathcal{R}} (\mathcal{U})$ as in Lemma \ref{lm:veroimgradealg}. 
Then the two natural morphisms 
$\mathbb{S} (\tilde{\mathcal{R}}(\mathcal{U})_1) \to 
\tilde{\mathcal{R}}(\mathcal{U})$ and 
$\mathbb{S} (\tilde{\mathcal{R}}(\mathcal{U})_2) \to 
\tilde{\mathcal{R}}(\mathcal{U})$ 
induce 
a rational map 
$\mathrm{Proj}\, \tilde{\mathcal{R}}(\mathcal{U}) 
-- \to \mathbb{P} (\mathcal{V}_1) |_{\mathcal{U}}$ 
and a closed immersion 
$\mathrm{Proj}\, \tilde{\mathcal{R}}(\mathcal{U}) 
\to \mathbb{P} (\mathcal{V}_2) |_{\mathcal{U}}$, respectively. 
Moreover, we see easily that the restriction 
$(\mathrm{Proj}\, \tilde{\mathcal{R}}(\mathcal{U})) 
|_{\mathcal{U} \setminus \mathrm{supp}\, \bar{\tau}}
\to 
\mathbb{P} (\mathcal{V}_1) 
|_{\mathcal{U} \setminus \mathrm{supp}\, \bar{\tau}}$
of this 
$\mathrm{Proj}\, \tilde{\mathcal{R}}(\mathcal{U}) 
-- \to \mathbb{P} (\mathcal{V}_1)
|_{\mathcal{U}}$
is an isomorphism. 
Thus, since $\mathrm{Proj}\, \tilde{\mathcal{R}}(\mathcal{U})$ 
is irreducible and reduced, we conclude that 
$\mathcal{W} |_{\mathcal{U}} = \mathrm{Proj}\, \tilde{\mathcal{R}}(\mathcal{U})$.
The morphism 
$\tilde{\mathcal{R}}(\mathcal{U})_4 \to 
\mathcal{R}(\mathcal{U})_4$ induced by the natural projection
$\tilde{\mathcal{R}}(\mathcal{U}) \to 
\mathcal{R}(\mathcal{U})$ coincides with the projection 
$\tilde{\mathcal{V}}_4 \to \mathcal{V}_4$ constructed above. \qed

In the description 
$\mathcal{W} |_{\mathcal{U}} = \mathrm{Proj}\, \tilde{\mathcal{R}}(\mathcal{U})$ 
of $\mathcal{W} |_{\mathcal{U}}$ above, the polynomial 
$y^2 - \bar{P}$ is a base of $\mathcal{L}_4^{\prime} (\mathcal{U})$ 
that cuts out $\mathcal{X} |_{\mathcal{U}}$ in $\mathcal{W} |_{\mathcal{U}}$. 
Meanwhile the natural rational map 
$\mathcal{W} |_{\mathcal{U}} -- \to \mathbb{P} (\mathcal{V}_1) |_{\mathcal{U}}$ 
is induced by the natural morphism 
$\mathbb{S} (\tilde{\mathcal{R}} (\mathcal{U})_1) 
\simeq \mathcal{O}_{\mathcal{B}} (\mathcal{U}) [x_0, x_1, x_2] 
\to \tilde{\mathcal{R}} (\mathcal{U})$, and so the 
polynomial $\bar{Q}^2 - h^2 \bar{P}$ is a base of 
$\mathcal{L}_4 (\mathcal{U})$ 
(i.e., the defining equation of $\varSigma |_{\mathcal{U}}$
in $\mathbb{P} (\mathcal{V}_1)$).   
Note that in $\tilde{\mathcal{R}} (\mathcal{U})$ we have

\begin{equation} \label{eql:basesl4l4prime}
\bar{Q}^2 - h^2 \bar{P} = h^2 (y^2 - \bar{P}). 
\end{equation}
Since $h$ is a defining equation of the divisor $\bar{\tau}$, 
it follows from (\ref{eql:basesl4l4prime}) that
\[
\mathcal{L}_4^{\prime} \simeq \mathcal{L}_4 
\otimes \mathcal{O}_{\mathcal{B}} (2 \bar{\tau}) .   
\]

\begin{lemma}   \label{lm:omegasbformula}
Let $\bar{\psi} : \mathcal{S} \to \mathcal{X}$
be the natural projection to the relative 
canonical model $\mathcal{X}$, 
$\bar{\psi}_1 : \mathcal{X} \to \varSigma$, 
the  natural projection to the relative $1$--canonical image
$\varSigma$, 
and 
$\frak{j}_{\varSigma} : \varSigma \to \mathbb{P} (\mathcal{V}_1)$,  
the natural inclusion to the projective bundle 
$\mathbb{P} (\mathcal{V}_1)$.
Let $\frak{g}: \mathcal{X} \to \mathcal{B}$, 
$\mathrm{pr}_{\varSigma} : \varSigma \to \mathcal{B}$, and 
$\mathrm{pr}_{\mathbb{P} (\mathcal{V}_1)} : 
\mathbb{P} (\mathcal{V}_1) \to \mathcal{B}$ be the 
natural projections. 
Let $\omega$ be the invertible sheaf on $\mathcal{X}$ defined by 
\[
 \omega = (\frak{j}_{\varSigma} \circ \bar{\psi_1})^*
\mathcal{O}_{\mathbb{P} (\mathcal{V}_1 )} 
(K_{\mathbb{P} (\mathcal{V}_1)}+ 
\varSigma - 
\mathrm{pr}_{\mathbb{P}(\mathcal{V}_1)}^* (K_{\mathcal{B}} + \bar{\tau}) ).
\]
Then $\omega_{\mathcal{S} | {\mathcal{B}}} \simeq \bar{\psi}^* \omega $, where
the sheaf $\omega_{\mathcal{S} | {\mathcal{B}}}$ is by definition 
$\mathcal{O}_{\mathcal{B}} (K_{\mathcal{S}} - \frak{f}^* K_{\mathcal{B}})$.  
\end{lemma}

Proof. 
Let $\mathcal{W}_1$ be the open set of $\mathcal{W}$ obtained by 
subtracting the indeterminacy locus of the rational map  
$\mathcal{W} - - \to \mathbb{P} (\mathcal{V}_1)$. 
Put $\mathcal{W}_0 = \mathcal{W}_1 - (\mathcal{W}_1)_{\mathrm{sing}}$. 
Let us denote by 
$\varrho : \mathcal{W}_1 \to \mathbb{P} (\mathcal{V}_1)$ 
the restriction to $\mathcal{W}_1$ of the rational map 
$\mathcal{W} - - \to \mathbb{P} (\mathcal{V}_1)$. 
Then as Cartier divisors the following hold:  
\begin{align}
\mathcal{X}|_{\mathcal{W}_1} 
&= (\varrho^* \varSigma )|_{\mathcal{W}_1} 
- 2 (\mathrm{pr}_{\mathcal{W}}^* (\bar{\tau}))|_{\mathcal{W}_1} 
\label{eql:xdivonw1} \\ 
K_{\mathcal{W}_0}     
&= (\varrho^* K_{\mathbb{P}(\mathcal{V}_1)})|_{\mathcal{W}_0}
+(\mathrm{pr}_{\mathcal{W}}^*(\bar{\tau}))|_{\mathcal{W}_0}, 
\label{eql:canow0} 
\end{align} 
where $\mathrm{pr}_{\mathcal{W}}: \mathcal{W} \to \mathcal{B}$ 
is the natural projection. 
Indeed (\ref{eql:xdivonw1}) follows from (\ref{eql:basesl4l4prime}), 
and (\ref{eql:canow0}) follows from a direct computation of 
pullbacks of differential forms using the descriptions 
in Lemma \ref{lm:veroimgradealg} of $\tilde{\mathcal{R}} (\mathcal{U})$. 
Adding (\ref{eql:xdivonw1}) and (\ref{eql:canow0}), we see that 
$(\frak{j}_{\varSigma} \circ \bar{\psi}_1)^* 
\mathcal{O}_{\mathbb{P}(\mathcal{V}_1)} 
(K_{\mathbb{P}(\mathcal{V}_1)} + \varSigma 
- \mathrm{pr}_{\mathbb{P}(\mathcal{V}_1)}^* (K_{\mathcal{B}}+\bar{\tau}))$
is an invertible sheaf on $\mathcal{X}$
whose restriction to $\mathcal{X}_0 
= \mathcal{X} \setminus \mathcal{X}_{\mathrm{sing}}$ 
coincides with the sheaf 
$\omega_{\mathcal{S} | \mathcal{B}} = 
\mathcal{O}_{\mathcal{S}} (K_{\mathcal{S}} - \frak{f}^* K_{\mathcal{B}} )$. 
Thus letting 
$\{ \mathcal{D}_i \}$ be the 
set of irreducible divisors on $\mathcal{S}$ exceptional with respect to 
$\bar{\psi} : \mathcal{S} \to \mathcal{X}$, we see that 
their exists a collection $\{ a_i \}$ ($a_i \in \mathbb{Z}$) such that 
$(\bar{\psi}^* \omega) \otimes (\omega_{\mathcal{S}|\mathcal{B}})^{\otimes (-1)}
\simeq \mathcal{O}_{\mathcal{S}} (\sum a_i \mathcal{D}_i)$ holds. 
Thus for any closed point $t \in T$ and for any 
$(-2)$--curve $D$ on $S_t$ exceptional with respect to 
$S_t \to X_t$, the equality 
$(\sum a_i \mathcal{D}_i|_{S_t}) D = 0$ holds.
Since the intersection form of a fundamental cycle of a rational 
point is nondegenerate, it follows from this that 
$a_i = 0$ for all $i$'s. Thus we obtain the assertion. \qed 
  
Now since 
$\varSigma \in |\mathcal{O}_{\mathbb{P} (\mathcal{V}_1)}(4) \otimes
\mathrm{pr}_{\mathbb{P} (\mathcal{V}_1)}^* (\mathcal{L}_4^{\otimes (-1)})|$,
we have 
\[
\omega \simeq 
(\frak{j}_{\varSigma} \circ \bar{\psi}_1)^*
(\mathcal{O}_{\mathbb{P} (\mathcal{V}_1)} (1)
\otimes 
\mathrm{pr}_{\mathbb{P} (\mathcal{V}_1)}^* 
((\det \mathcal{V}_1) \otimes 
\mathcal{O}_{\mathcal{B}} (- \bar{\tau}) \otimes \mathcal{L}_4^{\otimes (-1)})). 
\] 
From this it follows
\begin{equation}  \label{eql:vnformula}
\mathcal{V}_n \simeq \frak{g}_* (\omega^{\otimes n}) \simeq 
({\mathrm{pr}_{\varSigma}}_* 
{\bar{\psi_1}}_* \bar{\psi}_1^* \mathcal{O}_{\varSigma} (n))
\otimes
((\det \mathcal{V}_1) \otimes 
\mathcal{O}_{\mathcal{B}} (-\bar{\tau}) 
\otimes \mathcal{L}_4^{\otimes (-1)})^{\otimes n}, 
\end{equation}
where 
$\mathcal{V}_n \simeq \frak{g}_* (\omega^{\otimes n}) $ 
follows from Lemma \ref{lm:omegasbformula} 
and the normality of $\mathcal{X}$. 
Since we see by easy computation 
$({\mathrm{pr}_{\varSigma}}_* 
{\bar{\psi_1}}_* \bar{\psi}_1^* 
\mathcal{O}_{\varSigma} (n)) \simeq \mathcal{V}_n$, 
we obtain (by taking the determinant bundle of (\ref{eql:vnformula}) 
for $n=1$ and $n=4$) that  
\[
((\det \mathcal{V}_1) \otimes 
\mathcal{O}_{\mathcal{B}} (-\bar{\tau}) 
\otimes \mathcal{L}_4^{\otimes (-1)})^{\otimes 3} \simeq
((\det \mathcal{V}_1) \otimes 
\mathcal{O}_{\mathcal{B}} (-\bar{\tau}) 
\otimes \mathcal{L}_4^{\otimes (-1)})^{\otimes 56}
\simeq \mathcal{O}_{\mathcal{B}}, 
\]
hence $\mathcal{L}_4 \simeq (\det \mathcal{V}_1) \otimes 
\mathcal{O}_{\mathcal{B}} (-\bar{\tau})$. 
Thus we obtain the following: 

\begin{proposition}   \label{prop:l4prime}
$\mathcal{L}_4^{\prime} \simeq (\det \mathcal{V}_1) \otimes 
\mathcal{O}_{\mathcal{B}} (\bar{\tau})$. 
\end{proposition}

Now let us define the $5$--tuple for 
non-hyperelliptic deformation families.  
Let $T$ be a smooth variety,  
$t_0 \in T$, a closed point, and 
$\pi_{\mathcal{B}} : \mathcal{B} \to T$, 
a smooth proper family of irreducible curves of genus $\frak{b}$. 
Let $\mathcal{V}_1$ be a locally free sheaf of rank $3$, 
and $\bar{\tau}$, an effective divisor on $\mathcal{B}$ 
such that $B_{t_0} = \pi_{\mathcal{B}}^{-1} (t_0) 
\subset \mathrm{supp} \bar{\tau}$. 
Assume that we are given an invertible 
$\mathcal{O}_{\bar{\tau}}$--module $\bar{\mathcal{T}}$. 
Assume moreover that we are given an extension class
\[
\xi : \qquad 0 \to \mathbb{S}^2 (\mathcal{V}_1) \to \mathcal{V}_2
        \to \bar{\mathcal{T}} \to 0
\]
such that $\mathcal{V}_2$ is a locally free sheaf (of rank $6$) 
on $\mathcal{B}$. 
Then the morphism $\bar{\sigma_2} : 
\mathbb{S}^2 (\mathcal{V}_1) \to \mathcal{V}_2$ 
in the extension class $\xi$ induces 
a rational veronese map 
$\mathbb{P} (\mathcal{V}_1) - - \to \mathbb{P} (\mathcal{V}_2)$. 
We denote by $\mathcal{W}$ the image of this rational veronese map. 
Let $\bar{\frak{c}} : \mathbb{S}^2 (\bigwedge^2 \mathcal{V}_1) 
\to \mathbb{S}^2 (\mathbb{S}^2 (\mathcal{V}_1))$ be the morphism 
given in the paragraph just after Lemma \ref{lm:lnandsgmn}, 
and define the coherent sheaf 
$\tilde{\mathcal{V}_4}$ by 
$\tilde{\mathcal{V}_4} = 
\mathrm{Cok}\, (\mathbb{S}^2 (\bar{\sigma_2}) \circ \bar{\frak{c}})$. 
Then $\tilde{\mathcal{V}_4}$ is a locally free sheaf of rank $15$. 
Put $\mathcal{L}_4^{\prime}  = (\det \mathcal{V}_1) \otimes 
\mathcal{O}_{\mathcal{B}} (\bar{\tau})$, and assume that we are 
given an inclusion $\bar{\frak{j}} : 
\mathcal{L}_4^{\prime} \to \tilde{\mathcal{V}_4}$.
Then $\bar{\frak{j}}$ cuts out a hypersurface $\mathcal{X}$ 
of $\mathcal{W}$. Consider the following condition: 
\medskip

CD) For any closed point $b \in \mathrm{supp} \bar{\tau}$ 
the image of the morphism 
$\mathcal{V}_2 \otimes \mathbb{S}^2 (\mathcal{V}_1) \otimes k(b) \to  
\tilde{\mathcal{V}}_4 \otimes k(b)$ and the image of 
$\mathcal{L}_4^{\prime} \otimes k(b) 
\to \tilde{\mathcal{V}}_4 \otimes k(b)$  
intersect each other only at $0 \in \tilde{\mathcal{V}}_4 \otimes k(b)$, 
where 
$\mathcal{V}_2 \otimes \mathbb{S}^2 (\mathcal{V}_1)  \to  
\tilde{\mathcal{V}}_4$
is the composite of the morphism 
$(\mathrm{id}_{\mathcal{V}_2} \cdot \bar{\sigma}_2 ):
\mathcal{V}_2 \otimes \mathbb{S}^2 (\mathcal{V}_1)  
\to \mathbb{S}^2 (\mathcal{V}_2)$
and the natural projection
$\mathbb{S}^2 (\mathcal{V}_2) \to 
\tilde{\mathcal{V}}_4$. 
\medskip

We call $(\pi_{\mathcal{B}} : \mathcal{B} \to T,\,  
\mathcal{V}_1 ,\, \bar{\tau} ,\, \xi ,\, \bar{\frak{j}})$ a 
$5$--tuple for non--hyperelliptic deformation families, 
if it satisfies the condition CD) above. 
Moreover, we say that a $5$-tuple 
$(\pi_{\mathcal{B}} : \mathcal{B} \to T,\, 
\mathcal{V}_1 ,\, \bar{\tau} ,\, \xi ,\, \bar{\frak{j}})$ is 
admissible, if the following two conditions are satisfied: 
\medskip

I) let $\frak{g} : \mathcal{X} \to \mathcal{B}$ 
be the natural projection and put 
$\pi_{\mathcal{X}} =
\pi_{\mathcal{B}} \circ \frak{g}: \mathcal{X} \to T$; 
then for any closed point $t \in T$, 
the scheme theoretic fiber $X_t = \pi_{\mathcal{X}}^{-1} (t)$ 
has at most rational points as its singularities;  

II) $\frak{g} : \mathcal{X} \to \mathcal{B}$ admits a 
simultaneous resolution $\bar{\psi} : \mathcal{S} \to \mathcal{X}$ without 
taking a branch cover of $T$, such that for any closed point 
$t \in T$ the restriction ${\psi}_t :S_t \to X_t$ of $\bar{\psi}$ 
to $S_t$ gives the minimal resolution of $X_t$, where 
$\pi_{\mathcal{S}} = \pi_{\mathcal{X}} \circ \bar{\psi} 
: \mathcal{S} \to \mathcal{X}$ 
and $S_t = \pi_{\mathcal{S}}^{-1} (t)$. 
\medskip

\begin{remark}   \label{rm:trsversltycd}
The condition CD) is equivalent to say that 
$\mathcal{X}$ does not intersects the indeterminacy locus 
of $\mathcal{W} - - \to \mathbb{P} (\mathcal{V}_1)$.   
\end{remark}

Given an admissible $5$-tuple, let $\mathbb{S}^4 (\mathcal{V}_1)  
\to \tilde{\mathcal{V}}_4 / \mathcal{L}_4^{\prime}$ be 
the composite of the natural morphism 
$\mathbb{S}^4 (\mathcal{V}_1) 
\simeq \mathrm{ker}\, \bar{\frak{c}} \to 
\tilde{\mathcal{V}}_4 
\simeq \mathrm{ker}\, 
(\mathbb{S}^4 (\mathcal{V}_1) \circ \bar{\frak{c}})$ 
and the natural projection 
$\tilde{\mathcal{V}}_4 \to \tilde{\mathcal{V}}_4 / \mathcal{L}_4^{\prime}$. 
Define the invertible sheaf $\mathcal{L}_4$ by 
$\mathcal{L}_4 = \mathrm{ker}\, (\mathbb{S}^4 (\mathcal{V}_1)  
\to \tilde{\mathcal{V}}_4 / \mathcal{L}_4^{\prime})$. 
Then $\mathcal{L}_4$ cuts out a relative quartic  
$\varSigma \subset \mathbb{P} (\mathcal{V}_1)$.  
The natural projection $\bar{\psi}_1 : \mathcal{X} \to \varSigma$ 
is a birational morphism, and gives an isomorphism when restricted 
to the inverse image of $\mathcal{B} \setminus \mathrm{supp} \bar{\tau}$. 
By the same computation as in 
the proof of Lemma \ref{lm:omegasbformula}, 
we see that 
$\mathcal{O}_{\mathcal{S}} (K_{\mathcal{S}} - \frak{f}^* K_{\mathcal{B}} )
\simeq \bar{\psi}_1^* \mathcal{O}_{\varSigma} (1)$.  
From this it follows that 
$f_t : S_t \to B_t$ is a relatively minimal genus $3$ fibration. 
Since the image of $\bar{\sigma_2} \otimes \mathrm{id}_{k(b)}
: \mathbb{S}^2 (\mathcal{V}_1) \otimes k(b) \to 
\mathcal{V}_2 \otimes k(b)$ has rank $\geq 5$ for any closed 
point $b \in \mathcal{B}$, it follows that all the fibers 
of $f_t : S_t \to B_t$ are $2$--connected 
(see, e.g.,\cite[Lemma 4.1]{notesonI}).
We shall call 
$\frak{f} : \mathcal{S} \to \mathcal{B}$ 
the non-hyperelliptic deformation family 
(of the fibration $f_{t_0} : S_{t_0} \to B_{t_0}$) 
associated to the $5$--tuple 
$(\pi_{\mathcal{B}} : \mathcal{B} \to T,\,  
\mathcal{V}_1 ,\, \bar{\tau} ,\, \xi ,\, \bar{\frak{j}})$.

By Lemma \ref{lm:omegasbformula} and 
Proposition \ref{prop:l4prime}, we obtain the following: 

\begin{proposition}
Let $(\mathcal{S}, \mathcal{B}, T, t_0, 
\frak{f}, \pi_{\mathcal{S}},  \pi_{\mathcal{B}})$ 
be a $2$--connected non-hyperelliptic 
deformation family as in the beginning of this section 
(of a relatively minimal genus $3$ hyperelliptic fibration 
with all fibers $2$--connected). 
Put $\omega_{\mathcal{S} | \mathcal{B}} = 
\mathcal{O}_{\mathcal{S}} (K_{\mathcal{S}} - \frak{f}^* K_{\mathcal{B}})$ 
and $\mathcal{V}_n = 
\frak{f}_* (\omega_{\mathcal{S} | \mathcal{B}}^{\otimes n})$.   
Let $\bar{\sigma}_2 : \mathbb{S}^2 (\mathcal{V}_1) \to \mathcal{V}_2$
be the morphism induced by the multiplication structure of 
the relative canonical algebra of $\bar{\frak{f}}$. 
Let $\bar{\tau}$, $\bar{\mathcal{T}}$, and 
$\xi : 0 \to \mathbb{S}^2 (\mathcal{V}_1) \to \mathcal{V}_2 
\to \bar{\mathcal{T}} \to 0$ be as in the paragraph just before 
Lemma \ref{lm:rkimsigman}. 
Let moreover $\bar{\frak{j}} : \mathcal{L}_4^{\prime} \to 
\tilde{\mathcal{V}_4}$ be as in the paragraph just before 
Lemma \ref{lm:xhpsfofw}. 
Then $(\pi_{\mathcal{B}} : \mathcal{B} \to T,\,  
\mathcal{V}_1 ,\, \bar{\tau} ,\, \xi ,\, \bar{\frak{j}})$ is 
an admissible $5$--tuple.   
\end{proposition}  

By the computation above we obtain the following: 

\begin{theorem}   \label{thm:strtheorem}
Let $f : S \to B$ be a relatively minimal genus $3$ hyperelliptic 
fibration all of whose fibers are $2$--connected. Then 
the isomorphism classes of $2$--connected 
non-hyperelliptic deformation 
families (of the fibration $f : S \to B$) are in 
one--to--one correspondence with the isomorphism classes 
of admissible $5$--tuples for 
$2$--connected genus $3$ 
non-hyperelliptic 
deformation families for which 
the central fibrations $f_{t_0} : S_{t_0} \to B_{t_0}$'s coincide  
with $f : S \to B$. Here by a $2$--connected genus $3$ 
non-hyperelliptic deformation family, we mean a 
deformation family 
$(\mathcal{S}, \mathcal{B}, T, t_0, 
\frak{f}, \pi_{\mathcal{S}},  \pi_{\mathcal{B}})$ 
(in the sense of the beginning of this section) 
such that for any general closed point $t \in T$, the 
fibration $f_t: S_t = \mathrm{pr}_{\mathcal{S}}^{-1}(t) \to 
B_t = \mathrm{pr}_{\mathcal{B}}^{-1}(t)$ is non-hyperelliptic, 
and such that for any closed point $t \in T$, 
the fibration $f_t: S_t = \mathrm{pr}_{\mathcal{S}}^{-1}(t) \to 
B_t = \mathrm{pr}_{\mathcal{B}}^{-1}(t)$ is a relatively minimal 
genus $3$ fibration all of whose fibers are $2$--connected. 
\end{theorem}

\section{Some regular surfaces with $c_1^2 = 8$ and $p_g =4$}

\label{scn:applsfswith}

In this section, using Theorem \ref{thm:strtheorem}, we shall study 
some regular surfaces with $c_1^2 = 8$ and $p_g =4$. 
More precisely, 
in the moduli space of minimal regular surfaces with $c_1^2 = 8$ 
and $p_g =4$,  
we shall find two strata  
$\mathcal{M}_0^{\sharp}$ and $\mathcal{M}_0^{\flat}$
(each of dimensions $32$ and $30$, respectively)
whose members have no canonical involution, 
but have the same deformation type  
as that of members of the stratum $\mathcal{M}_0$ 
 in Bauer--Pignatelli \cite{caninvpg4c8}. 
It will turn out that 
the Stratum $\mathcal{M}_0$ 
(resp. a substratum of $\mathcal{M}_0$) is 
at the boundary of our stratum $\mathcal{M}_0^{\sharp}$ 
(resp. $\mathcal{M}_0^{\flat}$).  
\medskip

{\sc An existence theorem}

First, we prepare for an existence theorem for non-hyperelliptic 
deformation families of given $2$--connected genus $3$ 
hyperelliptic fibrations, using the language of $5$--tuples for 
($2$--connected) hyperelliptic genus $3$ fibrations given in 
the part I of this series \cite{notesonI}. For the definition 
of a $5$--tuple for hyperelliptic fibrations, see \cite{notesonI}. 

Assume that we are given a relatively minimal genus $3$ 
hyperelliptic fibration $f : S \to B$ with all fibers $2$-connected.
We put $V_n = f_* (\omega_{S | B}^{\otimes n})$ for each $n \geq 0$.
Then the hyperelliptic involution of $f$ induces a splitting 
$V_n = V_n^+ \oplus V_n^-$ into the eigen--sheaves $V_n^+$ and $V_n^-$. 
The sheaves $V_n^+$ and $V_n^-$ are locally free. 
We have $\mathrm{rk}\, V_1 = \mathrm{rk}\, V_1^- = 3$, 
$\mathrm{rk}\, V_n^{\pm} = 2n + 1$, and 
$\mathrm{rk}\, V_n^{\mp} = 2n - 3$, 
where $\pm$ stands for $+$ if $n$ is even, 
and for $-$ if $n$ is odd, 
and $\mp$ stands for $-$ if $n$ is even, 
and for $+$ if $n$ is odd 
(throughout this paper, 
we shall keep this rule as to what the symbols $\pm$ and $\mp$ mean).  
In particular we have 
$\mathrm{rk}\, V_2^+ = 5$ and $\mathrm{rk}\, V_2^- = 1$. 

Let $\sigma_n : \mathbb{S}^n (V_1) \to V_n$ be the 
natural morphism determined by the multiplication structure 
of the relative canonical algebra of 
$f: S \to B$,  and denote by 
$\sigma_n^{\pm} : \mathbb{S}^n (V_1)
= \mathbb{S}^n (V_1^-) \to V_n^{\pm}$ 
the restriction of $\sigma_n$ to $V_n^{\pm} \subset V_n$.  
Put $L = \mathrm{ker}\, \sigma_2^+$. 
For the proof of the following Lemma see \cite{notesonI}:
\begin{lemma}  \label{lm:v2minus}
$V_2^- \simeq (\det V_1) \otimes L^{\otimes (-1)}$. 
\end{lemma}

Now let us denote by 
$\frak{c} : \mathbb{S}^2 (\bigwedge^2 V_1) 
\to \mathbb{S}^2 (\mathbb{S}^2 (V_1))$ 
the morphism given by 
$(a \wedge b) (c \wedge d) \mapsto (ac)(bd) - (ad)(bc)$, 
and define the coherent sheaf $\tilde{V}_4$ on $B$ by 
\[
\tilde{V}_4 = 
\mathrm{Cok}\, (\mathbb{S}^2 (\sigma_2) \circ \frak{c} : 
\mathbb{S}^2 (\bigwedge^2 V_1) 
\to  
\mathbb{S}^2 (V_2)). 
\] 

Then the proof of Catanese--Pignatelli 
\cite[Lemma 7.7]{fibrationsI'} works 
also for our case of hyperelliptic fibrations. Thus we 
see that for any closed point $b \in B$ the morphism 
$\mathbb{S}^2 (\bigwedge^2 V_1) \otimes k(b) \to  
\mathbb{S}^2 (V_2) \otimes k(b) $ has rank $6$, 
and that the morphism 
$\mathbb{S}^2 (\sigma_2) \circ \frak{c}$ 
is injective. It follows that $\tilde{V}_4$ is a locally 
free sheaf of rank $\mathrm{rk}\, \tilde{V}_4 = 15$.  
Note that for our case, we have $V_2 = V_2^+ \oplus V_2^-$, 
and that the morphism $\mathbb{S}^2 (\sigma_2) \circ \frak{c}$ 
factors through the natural inclusion 
$\mathbb{S}^2 (V_2^+) \to \mathbb{S}^2 (V_2)$. 
Thus $\mathrm{Cok}\, (\mathbb{S}^2 (\sigma_2^+) \circ \frak{c})$ 
is a locally free sheaf on $B$ of rank $9$, and we have
\[
\tilde{V}_4 = 
\mathrm{Cok}\, (\mathbb{S}^2 (\sigma_2^+) \circ \frak{c}) \oplus
(V_2^+ \otimes V_2^-) \oplus 
(V_2^-)^{\otimes 2}. 
\]

Consider the natural morphism 
$\mathbb{S}^4 (V_1) \simeq \mathrm{Cok}\, \frak{c} \to 
\mathrm{Cok}\, (\mathbb{S}^2 (\sigma_2^+) \circ \frak{c})$. 
Since the composite 
$\mathbb{S}^2 (V_1) \otimes L \to 
\mathrm{Cok}\, (\mathbb{S}^2 (\sigma_2^+) \circ \frak{c})$
of the two natural morphisms  
$\mathbb{S}^2 (V_1) \otimes L \to \mathbb{S}^4 (V_1)$ and 
$\mathbb{S}^4 (V_1) \simeq \mathrm{Cok}\, \frak{c} \to 
\mathrm{Cok}\, (\mathbb{S}^2 (\sigma_2^+) \circ \frak{c})$ 
is a zero morphism, we obtain a natural projection 
\begin{equation} \label{eql:a4tocok}
A_4 \to \mathrm{Cok}\, (\mathbb{S}^2 (\sigma_2^+) \circ \frak{c}),
\end{equation}
where $A_4 = \mathrm{Cok}\, 
(\mathbb{S}^2 (V_1) \otimes L \to \mathbb{S}^4 (V_1))$. 
Since $A_4$ is also a locally free sheaf of rank $9$ 
(see \cite[Lemma 3.5]{notesonI}), this induced projection 
(\ref{eql:a4tocok}) is an isomorphism. 

Thus we obtain the following:
\begin{lemma}    \label{lm:v4isomoplus}
$\tilde{V}_4 \simeq 
A_4 \oplus (V_2^+ \otimes V_2^-) \oplus (V_2^-)^{\otimes 2}$.  
\end{lemma}

The next theorem gives a set of sufficient conditions for the 
existence of non--hyperelliptic deformation families. 
Note that for our study of regular surfaces with 
$c_1^2 = 8$ and $p_g= 4$, 
the explicit construction in the proof of the theorem    
is as important as the assertion itself of the theorem. 

\begin{theorem}   \label{thm:exstcthm}
Let $f: S \to B$ be a relatively minimal 
genus $3$ hyperelliptic fibration all of whose 
fibers are $2$--connected. 
Let $V_n = V_n^+ \oplus V_n^- $ 
be the decomposition of 
$f_* (\omega_{S | B}^{\otimes n}) = V_n$ 
into the eigen--sheaves with 
respect to the hyperelliptic involution.  
Put $L = \ker \sigma_2 \subset \mathbb{S}^2 (V_1)$, 
where $\sigma_2 : \mathbb{S}^2 (V_1) \to V_2$ is 
the multiplication morphism.
Put moreover  
$\tilde{V}_4 = \mathrm{Cok}\, (\mathbb{S}^2 (\sigma_2) \circ \frak{c})$
and $L_4^{\prime} = (V_2^-)^{\otimes 2}$. 
Assume  
1) $\dim \mathrm{Hom}\, (\mathbb{S}^2 (V_1),\, V_2^-) 
> \dim \mathrm{Hom}\, (V_2^+,\, V_2^-)$, 
2) $h^1 (\tilde{V}_4 \otimes (L_4^{\prime})^{\otimes (-1)}) = 0$, 
and 
3) that the relative canonical model $X$ of $S$ is non--singular.  
Then there exists a deformation family  
$(\mathcal{S}, \mathcal{B}, T, t_0, 
\frak{f}, \pi_{\mathcal{S}}, \pi_{\mathcal{B}})$ 
of the fibration $f_{t_0} = f : S_{t_0} = S \to B_{t_0} = B$ 
such that for any general closed point $t \in T$ 
the fibration $f_t : S_t = \pi_{\mathcal{S}}^{-1} (t) 
\to B_t = \pi_{\mathcal{B}}^{-1} (t)$ is non--hyperelliptic. 
\end{theorem} 

Proof. 
Take an affine neighbourhood $T \subset \mathbb{A}^1$ of 
$t_0 = \{ \zeta = 0  \} \in \mathbb{A}^1$, where 
$\zeta$ is the standard global coordinate of the affine line $\mathbb{A}^1$. 
In what follows, we frequently  replace the neighbourhood $T$ by  
smaller one without mentioning it explicitly.   

Put $\mathcal{B} = B \times T$, and denote by 
$\mathrm{pr}_1 : \mathcal{B} \to B$ 
and $\pi_{\mathcal{B}} : \mathcal{B} \to T$ 
the first projection and the second projection, respectively. 
We define the locally free sheaves $\mathcal{V}_1$, 
$\mathcal{V}_2^+$, and $\mathcal{V}_2^+$ on $\mathcal{B}$ 
by 
$\mathcal{V}_1 = \mathrm{pr}_1^* (V_1)$, 
$\mathcal{V}_2^+ = \mathrm{pr}_1^* (V_2^+)$, and 
$\mathcal{V}_2^- = \mathrm{pr}_1^* (V_2^-)$. 
Consider the natural short exact sequence 
\begin{equation}   \label{eql:shexcseqlsv2+}
0 \to L \to \mathbb{S}^2 (V_1) \to V_2^+ \to 0
\end{equation}
determined by the fibration $f: S \to B$.

By the exact sequence (\ref{eql:shexcseqlsv2+}) and 
the condition 1), there exists an element  
$\phi \in \mathrm{Hom} (\mathbb{S}^2 (V_1),\, V_2^- )$ 
such that 
$\phi_0 \neq 0 \in \mathrm{Hom} ( L,\, V_2^-)$, 
where $\phi_0 = \phi |_{L}$ is the restriction of 
$\phi$ to $L$. 
Using the two natural morphisms 
$\mathrm{pr}_1^*  \sigma_2^+ : 
\mathbb{S}^2 (\mathcal{V}_1) \to \mathcal{V}_2^+$ 
and 
$\mathrm{pr}_1^* \phi : 
\mathbb{S}^2 (\mathcal{V}_1) \to \mathcal{V}_2^-$ 
induced by  
$\sigma_2^+ : \mathbb{S}^2 (V_1) \to V_2^+$ and 
$\phi : \mathbb{S}^2 (V_1) \to V_2^-$, respectively, 
we define the morphism 
$\bar{\sigma}_2 : \mathbb{S}^2 (\mathcal{V}_1) \to 
\mathcal{V}_2 = \mathcal{V}_2^+ \oplus \mathcal{V}_2^-$
by $\bar{\sigma}_2 = 
(\mathrm{pr}_1^*  \sigma_2^+) \oplus ( \zeta \mathrm{pr}_1^* \phi )$. 
Note that for any non-zero constant $c \neq 0 \in \mathbb{C}$,  
the natural morphism 
$(\sigma_2^+ \oplus c \phi ) \otimes \mathrm{id}_{k(b)} : 
\mathbb{S}^2 (V_1) \otimes k(b) \to 
V_2 \otimes k(b)$ has 
rank $5$ if $b \in \mathrm{supp}\, \mathrm{Cok}\, (\phi_0)$, 
and has rank $6$ 
if $b \in B \setminus \mathrm{supp}\, \mathrm{Cok}\, (\phi_0)$.  
Thus there exist 
an effective divisor $\bar{\tau}$ on $\mathcal{B}$ 
and an invertible module $\mathcal{O}_{\bar{\tau}}$--module 
$\bar{\mathcal{T}}$ 
such that  
$B_{t_0} = \pi_{\mathcal{B}}^{-1} (t_0) 
\subset \mathrm{supp} \bar{\mathcal{T}}$ 
and 
$\mathrm{Cok}\, \bar{\sigma}_2 \simeq \bar{\mathcal{T}}$. 

Using this $\bar{\sigma}_2$, we define the locally free sheaf 
$\tilde{\mathcal{V}}_4$ on $\mathcal{B}$ by 
$\tilde{\mathcal{V}}_4 = 
\mathrm{Cok}\, (\mathbb{S}^2 (\bar{\sigma}_2) \circ \bar{\frak{c}} 
: \mathbb{S}^2 (\bigwedge^2 \mathcal{V}_1) \to 
\mathbb{S}^2 (\mathcal{V}_2)) $, where 
$\bar{\frak{c}} : \mathbb{S}^2 (\bigwedge^2 \mathcal{V}_1) \to
\mathbb{S}^2 (\mathbb{S}^2 (\mathcal{V}_1))$ is 
the morphism given in 
the paragraph just before Lemma \ref{lm:veroimgradealg}. 
Note that we have
\[
\tilde{\mathcal{V}}_4 |_{B_{t_0}}\simeq 
\mathrm{Cok}\, (\mathbb{S}^2 (\sigma_2) \circ \frak{c}) = 
\tilde{V}_4 = A_4 \oplus (V_2^+ \otimes V_2^-) \oplus (V_2^-)^{\otimes 2},  
\]
where $\tilde{V}_4$ is the sheaf on $B_{t_0}$ as in the beginning of 
this section.  
Let us also define the invertible sheaf $\mathcal{L}_4^{\prime}$ 
on $\mathcal{B}$ by 
$\mathcal{L}_4^{\prime} = (\det \mathcal{V}_1) \otimes 
\mathcal{O}_{\mathcal{B}} (\bar{\tau})$.   
Note that we have 
$\mathcal{L}_4^{\prime} |_{B_{t_0}} \simeq 
(\det V_1) \otimes \mathcal{O}_{B_{t_0}} (\tau)$, 
where $\tau$ is an effective divisor on $B= B_{t_0}$ 
such that $\mathcal{O}_{\tau} \simeq \mathrm{Cok}\, \phi_0 $.
Thus we have 
\[
\mathcal{L}_4^{\prime} |_{B_{t_0}} \simeq 
(\det V_1) \otimes V_2^- \otimes L^{\otimes -1} \simeq 
(V_2^-)^{\otimes 2}.
\]

Now let $\delta : (V_2^-)^{\otimes 2} \to A_4$ be the 
morphism determined by the multiplication structure 
of the relative canonical algebra of the fibration $f : S \to B$
 (see \cite[Remark 2]{notesonI}). 
We define the morphism of sheaves on $B$
\[
 \frak{j} : L_4^{\prime} \simeq (V_2^-)^{\otimes 2} \to 
\tilde{V}_4 = A_4 \oplus (V_2^+ \otimes V_2^-) \oplus 
(V_2^-)^{\otimes 2}, 
\]
by $\frak{j} = ( - \delta) \oplus 0 \oplus \mathrm{id}_{(V_2^-)^{\otimes 2}}$. 
This defines an element of 
$H^0 (\tilde{V}_4  \otimes (L_4^{\prime})^{\otimes (-1)}) 
= H^0 ((\tilde{\mathcal{V}}_4  \otimes 
(\mathcal{L}_4^{\prime})^{\otimes (-1)}) |_{B_{t_0}})$. 
Note that the sheaf $\tilde{\mathcal{V}}_4  \otimes 
(\mathcal{L}_4^{\prime})^{\otimes (-1)}$ is flat over $T$. 
Thus by the assumption 
$h^1 (\tilde{V}_4  \otimes (L_4^{\prime})^{\otimes (-1)}) = 0$ 
and the base change theorem, we obtain a morphism 
$\bar{\frak{j}} : \mathcal{L}_4^{\prime} \to \tilde{V}_4$ 
such that $\bar{\frak{j}} |_{B_{t_0}} = \frak{j}$.

Since 
$\frak{j} = (-\delta) \oplus 0 \oplus (\mathrm{id}_{(V_2^-)^{\otimes 2}})$ and 
$\bar{\frak{j}} |_{B_{t_0}} = \frak{j}$, we see that 
if we replace the affine neighbourhood $T$ by smaller one, 
the morphism  $\bar{\frak{j}}$ satisfies the condition CD) 
in Section \ref{scn:strthm} (see Remark \ref{rm:trsversltycd}).   
Thus $\pi_{\mathcal{B}} : \mathcal{B} \to T$, 
$\mathcal{V}_1$, $0 \to \mathbb{S}^2 (\mathcal{V}_1) 
\to \mathcal{V}_2 \to 
\bar{\mathcal{T}} \to 0$, $\bar{\tau}$, and $\bar{\frak{j}}$ 
obtained above form a $5$--tuple for non-hyperelliptic 
deformation families of the fibration $f_{t_0} = f: S_{t_0} = S 
\to B_{t_0} = B$. 

Now let $\frak{g}: \mathcal{X} \to \mathcal{B}$ be 
the relative canonical model associated to this $5$-tuple. 
By the definition of $\frak{j}$, we see that  
$\pi_{\mathcal{X}}^{-1} (t_0)$ 
(where $\pi_{\mathcal{X}} =  \pi_{\mathcal{B}}  \circ \frak{g}$) 
coincides with the relative canonical model $X$ of our $S$.   
Since we have assumed that $X$ is non--singular, 
we see that if we replace the affine neighbourhood $T$ by smaller one, 
our $\mathcal{X}$ is non--singular.
In particular, the $5$--tuple above is admissible. 
Thus, putting $\frak{f} = \frak{g}$ and $\mathcal{S} = \mathcal{X}$, 
we obtain the assertion.   \medskip   \qed  

\medskip

{\sc The stratum $\mathcal{M}_0^{\sharp}$}

Next let us construct a family of 
minimal regular surfaces with $c_1^2 = 8$ and $p_g = 4$ 
that allow non--hyperelliptic genus three fibrations. 
We use the same method as in Catanese--Pignatelli \cite{fibrationsI'}.
Let $B = \mathbb{P}^1$ be the projective line, and 
put $V_1 = 
\mathcal{O}_B (2) \oplus \mathcal{O}_B (2) \oplus \mathcal{O}_B (3)$  
and   
$V_2 = 
\mathcal{O}_B (5)^{\oplus 3} \oplus 
\mathcal{O}_B (5)^{\oplus 2} \oplus 
\mathcal{O}_B (6)$. 
Let $\{ \eta_0 , \eta_1 \}$ be a base of the vector space 
$H^0 (\mathcal{O}_B (1))$. 
Consider the natural inclusion
\begin{equation}   \label{eql:ntrincl}
\mathcal{O}_B (4)^{\oplus 3} \oplus 
\mathcal{O}_B (5)^{\oplus 2} \oplus 
\mathcal{O}_B (6)  
\subset  
V_2 =
\mathcal{O}_B (5)^{\oplus 3} \oplus 
\mathcal{O}_B (5)^{\oplus 2} \oplus 
\mathcal{O}_B (6)  
\end{equation} 
induced by the morphism 
$(\eta_0 ,\, \eta_0 - \eta_1,\, \eta_1 ) \times: 
\mathcal{O}_B (4)^{\oplus 3} \to \mathcal{O}_B (5)^{\oplus 3}$
and the identity of 
$\mathcal{O}_B (5)^{\oplus 2} \oplus 
\mathcal{O}_B (6)$. 
We take a general isomorphism
\begin{equation} \label{eql:s2v1toincl}
\mathbb{S}^2 (V_1) 
\simeq 
\mathcal{O}_B (4)^{\oplus 3} \oplus 
\mathcal{O}_B (5)^{\oplus 2} \oplus 
\mathcal{O}_B (6)  
\to 
\mathcal{O}_B (4)^{\oplus 3} \oplus 
\mathcal{O}_B (5)^{\oplus 2} \oplus 
\mathcal{O}_B (6) 
\end{equation}
and denote by $\sigma_2 : \mathbb{S}^2 (V_1) \to V_2$
the composite of the two morphisms 
(\ref{eql:s2v1toincl}) and (\ref{eql:ntrincl}). 
Then $\eta_0 \eta_1 (\eta_0 - \eta_1)= 0$ defines 
an effective divisor $\tau$ (of degree $3$) 
of $B$ such that 
$0 \to \mathbb{S}^2 (V_1) \to V_2 
\to \mathcal{T} = \mathcal{O}_{\tau} \to 0$ 
is exact.

Let $W^{\sharp}$ be the image of the rational veronese map 
$\mathbb{P} (V_1) - - \to \mathbb{P} (V_2)$ induced by $\sigma_2$. 
Then by the same method as in \cite[Theorem 8.2]{fibrationsI'}, 
we see that if the isomorphism (\ref{eql:s2v1toincl}) is 
sufficiently general, then the variety $W^{\sharp}$ has 
at most singularities supported on the indeterminacy 
locus of the natural birational map 
$W^{\sharp} - - \to \mathbb{P} (V_1)$. 
Let us define the invertible sheaf $L_4^{\prime}$ on $B$ by 
$L_4^{\prime} \simeq (\det V_1) \otimes \mathcal{O}_B ( \tau ) 
\simeq \mathcal{O}_B (10)$. 
Letting $\frak{c} : \mathbb{S}^2 (\bigwedge^2 V_1) \to 
\mathbb{S}^2 (\mathbb{S}^2 (V_1))$  
to be the morphism given by 
$(a \wedge b) (c \wedge d) \to (ac)(bd) - (ad)(bc)$, 
we again define the sheaf $\tilde{V}_4$ on $B$ by 
$\tilde{V}_4 = \mathrm{Cok}\, (\mathbb{S}^2 (\sigma_2) \circ \frak{c} : 
\mathbb{S}^2 (\mathbb{S}^2 (V_1)) \to \mathbb{S}^2 (V_2))$.  

\begin{lemma}
Let $\mathrm{pr}_{\mathbb{P}(V_2)} : \mathbb{P}(V_2) \to B$ 
be the structure morphism of $\mathbb{P}(V_2)$, 
and $L_4^{\prime} \to \tilde{V}_4$,  
a sufficiently general morphism.  
Then the hypersurface 
$X^{\sharp}$ of $W^{\sharp}$ cut out by the corresponding 
section of  
$H^0 ((\mathcal{O}_{\mathbb{P}(V_2)} (2) \otimes 
\mathrm{pr}_{\mathbb{P}(V_2)}^* ((L_4^{\prime})^{\otimes (-1)})$
$)|_{W^{\sharp}} )$
is a minimal regular surface with 
$c_1^2 = 8$ and $p_g = 4$. 
The natural projection $f^{\sharp} : S^{\sharp} = X^{\sharp} \to B$ 
gives a non--hyperelliptic genus $3$ fibration all of whose 
fibers are $2$--connected. 
\end{lemma}

Proof. 
Since we have 
$\mathbb{S}^2 (V_2) \simeq 
\mathcal{O}_B (10)^{\oplus 15} \oplus 
\mathcal{O}_B (11)^{\oplus 5} \oplus \mathcal{O}_B (12)$ 
and $L_4^{\prime} \simeq \mathcal{O}_B (10)$, 
the linear system 
$|(\mathcal{O}_{\mathbb{P}(V_2)} (2) \otimes 
\mathrm{pr}_{\mathbb{P}(V_2)}^* ((L_4^{\prime})^{\otimes (-1)}))|_{W^{\sharp}}|$  
is free from base points. 
Thus a general member $S^{\sharp} = X^{\sharp}$ of this linear 
system on $W^{\sharp}$ is non--singular. 
Moreover, since $f^{\sharp}_* (\omega_{S^{\sharp}}) 
\simeq \mathcal{O}_B \oplus \mathcal{O}_B \oplus 
\mathcal{O}_B (1)$, the canonical system 
$|\omega_{S^{\sharp}}|$ is free from base points. 
Thus the assertion follows from \cite[Theorem 7.13]{fibrationsI'}. 
\qed 

Let 
$\varPhi_{|\mathcal{O}_{\mathbb{P} (V_1)} (1) 
\otimes \mathrm{pr}_{\mathbb{P} (V_1)}^* ( \omega_B )|} 
: \mathbb{P} (V_1) \to \mathbb{P}^3$ 
be the morphism determined by the linear system 
$|\mathcal{O}_{\mathbb{P} (V_1)} (1) 
\otimes \mathrm{pr}_{\mathbb{P} (V_1)}^* ( \omega_B )|$. 
Let $\varSigma^{\sharp} \subset \mathbb{P} (V_1)$ be  
the relative $1$--canonical image of $S^{\sharp}$. 
Letting $x_0$, $x_1$, and $x_2$ be local bases of 
the direct summands 
$\mathcal{O}_B (2)$, $\mathcal{O}_B (2)$, and $\mathcal{O}_B (3)$   
of $V_1$, respectively,  
we see that 
$\varPhi_{|\mathcal{O}_{\mathbb{P} (V_1)} (1) 
\otimes \mathrm{pr}_{\mathbb{P} (V_1)}^* ( \omega_B )|}$ 
is the contraction of 
$\{ x_2 = 0 \} \simeq 
\mathbb{P} (\mathcal{O}_B (2) \oplus \mathcal{O}_B (2)) 
\subset \mathbb{P} (V_1)$ to a line in $\mathbb{P}^3$, 
and that $\varSigma^{\sharp}$ is not contained in 
$\{ x_2 = 0 \}$. 
Thus, since the canonical map $\varPhi_{|K_{S^{\sharp}}|}$ of 
$S^{\sharp}$ is the composite of the 
natural projection 
$S^{\sharp} = X^{\sharp} \to \varSigma^{\sharp} \subset \mathbb{P} (V_1)$ 
and the birational morphism 
$\varPhi_{|\mathcal{O}_{\mathbb{P} (V_1)} (1) 
\otimes \mathrm{pr}_{\mathbb{P} (V_1)}^* ( \omega_B )|}$, 
we see that the canonical map 
$\varPhi_{|K_{S^{\sharp}}|}$ is birational. 

Now let us compute the dimension of the stratum  
(in the moduli space) given by the surfaces 
with $c_1^2 = 8$ and $p_g = 4$ constructed above.  


First, note that for general $\sigma_2$ and 
general $L_4^{\prime} \to \tilde{V}_4$, 
the singular locus of $\varSigma^{\sharp}$ is 
$\varSigma^{\sharp} \cap 
\mathrm{pr}_{\mathbb{P}(V_1)}^{-1} (\mathrm{supp}\, \tau)$, 
and that for each $b \in \mathrm{supp}\, \tau$ the 
inverse image 
$\varSigma^{\sharp} \cap 
\mathrm{pr}_{\mathbb{P}(V_1)}^* (b)$ 
is a double conic in 
$\mathrm{pr}_{\mathbb{P}(V_1)}^{-1} (b) 
\simeq \mathbb{P}^2$ by (\ref{eql:basesl4l4prime}). 
It follows from this that 
the image 
$\varPhi_{|\mathcal{O}_{\mathbb{P}(V_1)} (1) 
\otimes \mathrm{pr}_{\mathbb{P}(V_1)}^* (\omega_B)|} 
(\{ x_2 = 0 \})$ is the unique 
singular line (in fact a quadruple line) of 
the canonical image 
$\varPhi_{|K_S^{\sharp}|} (S^{\sharp})$. 
The genus $3$ fibration $f^{\sharp}: S^{\sharp} \to B$ 
is given by the pullback of the linear system 
consisting of all the hyperplanes in $\mathbb{P}^3$ containing this line. 
Thus in order to compute the dimension of the stratum in the moduli space, 
we only need to compute the dimension of the family of the 
isomorphism classes (up to automorphisms of the base curve $B$) of 
the non--hyperelliptic genus $3$ fibrations constructed above.

Let us compute this dimension.  
%
%
%
For this, we count the number of parameters that we use in the 
construction.

First we may assume that 
$\mathrm{supp} \, \tau = \{ 0, 1, \infty \} \subset B = \mathbb{P}^1$.
Then the choice of the extension class of 
$0 \to \mathbb{S}^2 (V_1) \to V_2 \to \mathcal{O}_{\tau} \to 0 $ 
depends on 
$\dim (\mathrm{Ext}_B^1 (\mathcal{O}_{\tau}, \mathbb{S}^2 (V_1)) 
/ \mathrm{Aut} (\mathcal{O}_{\tau}))$ parameters. 
Applying the contravariant functor 
$\{ \mathrm{Ext}_B^n (\, \cdot \, , \mathbb{S}^2 (V_1)) \}_n$ to 
the short exact sequence 
$0 \to \mathcal{O}_B (-3)  \to \mathcal{O}_B \to \mathcal{O}_{\tau} \to 0$, 
we obtain 
$\dim \mathrm{Ext}_B^1 (\mathcal{O}_{\tau}, \mathbb{S}^2 (V_1)) = 18$.
Thus the choice of the extension class 
depends on $15$ parameters. 

Next, under a fixed extension class,  
the choice of $X^{\sharp}$ cut out by 
$L_4^{\prime} \to \tilde{V}_4$ depends on 
$\dim \mathbb{P} (H^0 (\tilde{V}_4 \otimes (L_4^{\prime})^{\otimes (-1)} ))$ 
parameters. Let us compute this dimension. 
Note that we have the short exact sequence 
$0 \to \mathbb{S}^2 (\bigwedge^2 V_1)  
  \to \mathbb{S}^2 (V_2) 
  \to \tilde{V}_4 
  \to 0$ 
tensored by $(L_4^{\prime})^{\otimes (-1)}$. 
Thus when 
$V_2 \simeq 
\mathcal{O}_B (5)^{\oplus 3} \oplus 
\mathcal{O}_B (5)^{\oplus 2} \oplus 
\mathcal{O}_B (6)$, 
we obtain 
$h^1 (\tilde{V}_4 \otimes (L_4^{\prime})^{\otimes (-1)}) = 0$. 
This implies 
$h^0 (\tilde{V}_4 \otimes (L_4^{\prime})^{\otimes (-1)})
= 
\chi (\mathbb{S}^2 (V_2) \otimes (L_4^{\prime})^{\otimes (-1)})
- \chi (\mathbb{S}^2 (\bigwedge^2 V_1) \otimes (L_4^{\prime})^{\otimes (-1)})$
and $h^1 (\tilde{V}_4 \otimes (L_4^{\prime})^{\otimes (-1)}) = 0$ 
for a general $S^{\sharp}$ in our construction.  
But we see easily that 
$\deg (\mathbb{S}^2 (V_2) \otimes (L_4^{\prime})^{\otimes (-1)}) = 7$ 
and 
$\deg (\mathbb{S}^2 (\bigwedge^2 V_1) \otimes (L_4^{\prime})^{\otimes (-1)})
= -4$. 
It follows that 
$h^0 (\tilde{V}_4 \otimes (L_4^{\prime})^{\otimes (-1)}) = 26$
for general $S^{\sharp}$. 
Thus the choice of $X^{\sharp}$ depends on $25$ parameters. 

So we have $15 + 25 = 40$ parameters. 
We however need to subtract 
$\dim \mathrm{Aut} (\mathbb{P} (V_1) / B) = 8$. 
From this we see that 
the minimal surfaces $S^{\sharp}$'s constructed above 
form a $32$--dimensional stratum 
in the moduli space of minimal regular surfaces with 
$c_1^2 = 8$ and $p_g = 4$.  
Thus we obtain the following:

\begin{proposition}      \label{prop:m0sharp}
In the moduli space of minimal regular surfaces with 
$c_1^2 = 8$ and $p_g = 4$, there exists a $32$--dimensional irreducible 
stratum $\mathcal{M}_0^{\sharp}$ whose members $S^{\sharp}$'s 
admit non--hyperelliptic genus $3$ fibrations 
$f^{\sharp} : S^{\sharp} \to B$'s with 
$V_1 \simeq 
\mathcal{O}_B (2) \oplus 
\mathcal{O}_B (2) \oplus 
\mathcal{O}_B (3)$ 
and 
$L_4^{\prime} \simeq \mathcal{O}_B (10)$, 
and such that the members $S^{\sharp}$'s have birational canonical maps. 
\end{proposition}

\begin{remark}       \label{rm:refmsharp}
It is likely that the canonical image of a general $S^{\sharp}$ 
constructed above belongs to Ciliberto's family $\mathscr{K} (8)$ 
in \cite[Theorem 8.1]{cilcansfpg4}.
We however remark that our $S^{\sharp}$ has canonical image 
with non-ordinary singularities (indeed, the quadruple line 
at the image by 
$\varPhi_{|\mathcal{O}_{\mathbb{P}(V_1)} (1) 
\otimes \mathrm{pr}_{\mathbb{P}(V_1)}^* (\omega_B)|}$
of  $\{ x_2 = 0 \}$), 
hence at least not corresponding  
to a main member of $\mathscr{K} (8)$. 
\end{remark}


\medskip

{\sc The stratum $\mathcal{M}_0^{\flat}$}

Let us construct one more family of minimal regular surfaces 
with $c_1^2 = 8$ and $p_g = 4$ that allow non--hyperelliptic 
genus $3$ fibrations. 
Let $B = \mathbb{P}^1$ again be the projective line, 
and put this time 
$V_1 = 
\mathcal{O}_B (1) \oplus \mathcal{O}_B (3) \oplus \mathcal{O}_B (3)$  
and   
$V_2 = 
\mathcal{O}_B (5) \oplus 
\mathcal{O}_B (4)^{\oplus 2} \oplus 
\mathcal{O}_B (6)^{\oplus 3}$. 
Consider the natural inclusion
\begin{equation}   \label{eql:ntrinclagain}
\mathcal{O}_B (2) \oplus 
\mathcal{O}_B (4)^{\oplus 2} \oplus 
\mathcal{O}_B (6)^{\oplus 3}  
\subset  
V_2 =
\mathcal{O}_B (5) \oplus 
\mathcal{O}_B (4)^{\oplus 2} \oplus 
\mathcal{O}_B (6)^{\oplus 3}  
\end{equation} 
induced by the morphism 
$\eta_0 \eta_1 ( \eta_0 - \eta_1  ) \times: 
\mathcal{O}_B (2) \to \mathcal{O}_B (5)$
and the identity of 
$\mathcal{O}_B (4)^{\oplus 2} \oplus 
\mathcal{O}_B (6)^{\oplus 3}$, 
where $\{ \eta_0 , \eta_1 \}$ is again a base of the 
vector space $H^0 (\mathcal{O}_B (1))$.
We take a general isomorphism
\begin{equation} \label{eql:s2v1toinclagain}
\mathbb{S}^2 (V_1) 
\simeq 
\mathcal{O}_B (2) \oplus 
\mathcal{O}_B (4)^{\oplus 2} \oplus 
\mathcal{O}_B (6)^{\oplus 3}  
\to 
\mathcal{O}_B (2) \oplus 
\mathcal{O}_B (4)^{\oplus 2} \oplus 
\mathcal{O}_B (6)^{\oplus 3}, 
\end{equation}
and repeat the same steps as in the case of 
the stratum $\mathcal{M}_0^{\sharp}$. 
Then the image $W^{\flat}$ of the rational veronese 
map $\mathbb{P}(V_1) - - \to \mathbb{P} (V_2)$ 
has at most singularities supported on the 
indeterminacy locus of the rational map 
$W^{\flat} - - \to \mathbb{P} (V_1)$. 

This time, the linear system 
$|\mathcal{O}_{\mathbb{P}(V_2)} (2) \otimes 
\mathrm{pr}_{\mathbb{P}(V_2)}^* ((L_4^{\prime})^{\otimes (-1)})|$
has non--empty base locus. 
But by the same method as in \cite[Theorem 8.2]{fibrationsI'}, 
we can show that the linear system 
$|\mathcal{O}_{\mathbb{P}(V_2)} (2) \otimes 
\mathrm{pr}_{\mathbb{P}(V_2)}^* ((L_4^{\prime})^{\otimes (-1)})||_{W^{\flat}}$
(and hence also
$|(\mathcal{O}_{\mathbb{P}(V_2)} (2) \otimes 
\mathrm{pr}_{\mathbb{P}(V_2)}^* ((L_4^{\prime})^{\otimes (-1)}))|_{W^{\flat}}|$)
is free from base points. 
Thus, putting 
$S^{\flat} = X^{\flat}$, we obtain a relatively minimal 
non--hyperelliptic genus $3$ fibration 
$f^{\flat} : S^{\flat} \to B$ all of whose fibers are $2$-connected, 
where 
$X^{\flat}$ is the relative canonical model associated to 
the $5$--tuple.  
  
Let $x_0$, $x_1$, and $x_2$ be local bases of the direct summands
$\mathcal{O}_B (1)$, $\mathcal{O}_B (3)$, and $\mathcal{O}_B (3)$, 
respectively, of $V_1$. 
Then the section 
$\varDelta = \{ x_1 = x_2 = 0 \} \simeq \mathbb{P}^1 \subset 
\mathbb{P} (V_1)$ of the projection 
$\mathrm{pr}_{\mathbb{P}(V_1)} : \mathbb{P} (V_1) \to B$ 
is the base locus of the linear system 
$|\mathcal{O}_{\mathbb{P}(V_1)} (1) \otimes 
\mathrm{pr}_{\mathbb{P}(V_1)}^* \omega_B|$. 
Note that we have a natural isomorphism 
$H^0 (\mathcal{O}_{\mathbb{P}(V_1)} (1) \otimes 
\mathrm{pr}_{\mathbb{P}(V_1)}^* \omega_B) \simeq 
H^0 (\omega_{S})$.  
Note also that if we let 
$Z_0$ be a local base of the direct summand 
$\mathcal{O}_B (5)$ of the sheaf 
$V_2 =
\mathcal{O}_B (5) \oplus 
\mathcal{O}_B (4)^{\oplus 2} \oplus 
\mathcal{O}_B (6)^{\oplus 3}$, 
then the coefficient to $Z_0^2$ of a 
general element of 
$H^0 (\mathbb{S}^2 (V_2) \otimes (L_4^{\prime})^{\otimes (-1)})$ 
is nowhere vanishing (because $L_4^{\prime} \simeq \mathcal{O}_B (10)$). 
From these we see easily that 
for general $S^{\flat}$ in our construction, 
the relative $1$--canonical image 
$\varSigma^{\flat} \subset \mathbb{P} (V_1)$ does not intersects 
the section $\varDelta \subset \mathbb{P} (V_1)$. 
In particular, we see that the canonical system $|K_{S^{\flat}}|$  
is free from base points, and that our 
$S^{\flat}$ is a minimal regular surface with 
$c_1^2 = 8$ and $p_g = 4$.  
Thus, we obtain the following: 

\begin{lemma} 
Let $\mathrm{pr}_{\mathbb{P}(V_2)} : \mathbb{P}(V_2) \to B$ 
be the structure morphism of $\mathbb{P}(V_2)$, 
and $L_4^{\prime} \to \tilde{V}_4$,   
a sufficiently general morphism.  
Then the hypersurface 
$X^{\flat}$ of $W^{\flat}$ cut out by the corresponding 
section of 
$H^0 ((\mathcal{O}_{\mathbb{P}(V_2)} (2) \otimes 
\mathrm{pr}_{\mathbb{P}(V_2)}^* ((L_4^{\prime})^{\otimes (-1)})$
$)|_{W^{\flat}} )$
is a minimal regular surface with 
$c_1^2 = 8$ and $p_g = 4$. 
The natural projection $f^{\flat} : S^{\flat} = X^{\flat} \to B$ 
gives a non--hyperelliptic genus $3$ fibration all of whose 
fibers are $2$--connected. 
\end{lemma}

Let 
$\varPhi_{|\mathcal{O}_{\mathbb{P} (V_1)} (1) 
\otimes \mathrm{pr}_{\mathbb{P} (V_1)}^* ( \omega_B )|} 
: \mathbb{P} (V_1) - - \to \mathbb{P}^3$ 
be the rational map determined by the linear system 
$|\mathcal{O}_{\mathbb{P} (V_1)} (1) 
\otimes \mathrm{pr}_{\mathbb{P} (V_1)}^* ( \omega_B )|$. 
This time we see easily that 
its image is a non--singular quadric 
$\mathcal{Q} \simeq \mathbb{P}^1 \times B \subset \mathbb{P}^3$, 
that this $\mathbb{P} (V_1) - - \to \mathcal{Q} \simeq \mathbb{P}^1 \times B$ 
is compatible with the projections to $B$, 
and that for each closed point $b \in B$ the restriction 
$\mathbb{P}^2 \simeq \mathrm{pr}_{\mathbb{P} (V_1)}^{-1} (b) 
- - \to \mathbb{P}^1 \times \{ b\} \subset \mathbb{P}^1 \times B$ 
is the linear projection 
from the point $\varDelta \cap \mathrm{pr}_{\mathbb{P} (V_1)}^{-1} (b)
\in \mathrm{pr}_{\mathbb{P} (V_1)}^{-1} (b)$. 
Since the canonical map $\varPhi_{|K_{S^{\flat}}|}$ of 
$S^{\flat}$ is the composite of 
the natural projection 
$S^{\flat} \to \varSigma^{\flat} \subset \mathbb{P} (V_1)$  
and the rational map 
$\varPhi_{|\mathcal{O}_{\mathbb{P} (V_1)} (1) 
\otimes \mathrm{pr}_{\mathbb{P} (V_1)}^* ( \omega_B )|} $,
it follows that 
$\varPhi_{|K_{S^{\flat}}|}$ is a  morphism of degree $4$ onto 
its image $\mathcal{Q}$. 

\begin{lemma}
If the surface $S^{\flat}$ in the construction above is general, 
then $S^{\flat}$ has no canonical involution. 
\end{lemma}

Proof. 
First let us prove that any general $S^{\flat}$ in our 
construction satisfies the following: 
for a general fiber $F$ of $f^{\flat} : S^{\flat} \to B$, 
the restriction 
$F \to \mathbb{P}^1 \times \{ f^{\flat}(F) \}
\subset \mathbb{P}^1 \times B$  
of $\varPhi_{|K_{S^{\flat}}|}$ has exactly $12$ critical 
values, and for each of these critical values there 
exists only one critical point of 
$F \to \mathbb{P}^1 \times \{ f^{\flat}(F) \}$ lying over it. 

For this purpose, let $\sigma_2 : \mathbb{S}^2 (V_1) \to V_2$ 
be the composite of 
the isomorphism (\ref{eql:s2v1toinclagain}) and 
the natural inclusion (\ref{eql:ntrinclagain}).  
Then, by $L_4^{\prime} = (\det V_1) \otimes \mathcal{O}_B (\tau) 
\simeq \mathcal{O}_B (10)$,  
we see easily that there exist  
an isomorphism (\ref{eql:s2v1toinclagain}) 
and a morphism 
$L_4^{\prime} \to \mathbb{S}^2 (V_2)$ ($\to \tilde{V}_4$) 
such that the associated $f^{\flat} : S^{\flat} \to B$ 
has a general fiber of the form 
$\{ x_0^4 + \alpha x_0 + \beta = 0  \} \subset 
\mathrm{pr}_{\mathbb{P} (V_1)}^{-1} (f^{\flat} (F)) \simeq \mathbb{P}^2$,   
where $\alpha \in k [x_1, x_2]_3$ and $\beta \in k [x_1, x_2]_4$ 
are general homogeneous polynomials of degree $3$ and $4$, 
respectively. 
Since the discriminant of the degree $4$ polynomial 
$x_0^4 + \alpha x_0 + \beta$ in $x_0$ is 
$u ((16)^2 \beta^3 - 27 \alpha^4)$ for a certain non--zero 
constant $u \in \mathbb{C}$, the restriction 
$\varPhi_{|K_{S^{\flat}}|} |_{F} : (x_0 : x_1 : x_2) \mapsto 
(x_1 : x_2) \in \mathbb{P}^1 \times \{ f^{\flat}(F) \} \subset 
\mathcal{Q}$ of the canonical map $\varPhi_{|K_{S^{\flat}}|}$ to 
$F$ has $12$ distinct critical values. 
This together with the Hurwitz formula implies that, also for general 
$\sigma_2 : \mathbb{S}^2 (V_1) \to V_2$ (with $V_1 \simeq 
\mathcal{O}_B (1) \oplus \mathcal{O}_B (3) \oplus \mathcal{O}_B (3)$), 
the associated 
$f^{\flat} : S^{\flat} \to B$ satisfies the condition stated at 
the beginning of this proof.      

Now assume that a general $S^{\flat}$ in our construction 
has a canonical involution $\iota$. 
Then, since the canonical map 
$\varPhi_{|K_{S^{\flat}}|} : S^{\flat} \to \mathcal{Q} \subset \mathbb{P}^3$
is the composite of the natural projection $S^{\flat} \to \varSigma^{\flat}$ 
and the rational map $\varPhi_{|\mathcal{O}_{\mathbb{P} (V_1)} (1) 
\otimes \mathrm{pr}_{\mathbb{P} (V_1)}^* ( \omega_B )|} : \mathbb{P} (V_1) 
- - \to \mathbb{P}^3$, any general fiber $F$ of $f^{\flat}$ is 
stable under the action by $\langle \iota \rangle \simeq 
\mathbb{Z} / 2$, where $\langle \iota \rangle$ is the group  
generated by the involution $\iota$. 
The natural projection $F \to F / \langle \iota \rangle$ is 
a double cover onto a non--singular curve 
$F / \langle \iota \rangle$ 
(of genus $g (F / \langle \iota \rangle)$). 
Recall that the restriction 
$\varPhi_{|K_{S^{\flat}}|} |_F : F \to \varPhi_{|K_{S^{\flat}}|} (F)$ 
to $F$ of the canonical map $\varPhi_{|K_{S^{\flat}}|}$ has 
$12$ critical values. By the uniqueness of the critical point 
lying over each of these $12$ critical values, we see that 
all of these $12$ critical points are fixed points of $\iota$, 
which contradicts the implication by the Hurwitz formula that  
the number of the fixed points of $\iota$ is 
$8 - 4 g (F / \langle \iota \rangle)$.   
Thus our general $S^{\flat}$ has no canonical involution. \qed

Now let us compute the dimension of the stratum (in the 
moduli space of minimal regular surfaces with 
$c_1^2 = 8$ and $p_g = 4$) filled up by $S^{\flat}$'s of 
our construction. 
We use the same method as in the computation for the 
case of the stratum $\mathcal{M}_0^{\sharp}$. 


First, recall that for general $\sigma_2$ and general 
$L_4^{\prime} \to \tilde{V}_4$, the canonical image 
$\varPhi_{|K^{S^{\sharp}}|} (S^{\sharp})$ of $S^{\sharp}$ is 
a non--singular quadric  
$\mathcal{Q} \simeq \mathbb{P}^1 \times B \subset \mathbb{P}^3$. 
The fibration $f^{\flat} : S^{\flat} \to B$ is given by 
the pullback of one of the two rulings of the quadric 
$\mathcal{Q}$. 
Thus in order to compute the dimension of the stratum, 
we only need to compute the dimension of the family 
of isomorphism classes (up to automorphisms of the 
base curve $B$) of the non--hyperelliptic genus $3$ fibrations 
constructed above.  

Let us compute this dimension.

%

First, by the same argument as in the case 
of the stratum $\mathcal{M}_0^{\sharp}$, 
we obtain 
$\dim \mathrm{Ext}_B^1 (\mathcal{O}_{\tau}, \mathbb{S}^2 (V_1)) = 18$. 
Thus the choice of the extension class 
$0 \to \mathbb{S}^2 (V_1) \to V_2 \to \mathcal{O}_{\tau} \to 0$ 
depends on $15$ parameters. 

Next, under the fixed extension class, the choice of the 
hypersurface $X^{\flat}$ of $W^{\flat}$ cut out by 
$L_4^{\prime} \to \tilde{V}_4$ depends on 
$h^0 (\tilde{V}_4 \otimes (L_4^{\prime})^{\otimes (-1)})$ 
parameters. Let us compute this dimension. 
Note that we have 
$V_1 \simeq 
\mathcal{O}_B (1) \oplus \mathcal{O}_B (3) \oplus \mathcal{O}_B (3)$, 
hence $\bigwedge^2 V_1 \simeq 
\mathcal{O}_B (4)^{\oplus 2} \oplus \mathcal{O}_B (6)$ 
and 
$\mathbb{S}^2 (\bigwedge^2 V_1) \simeq 
\mathcal{O}_B (8)^{\oplus 3} \oplus
\mathcal{O}_B (10)^{\oplus 2} \oplus
\mathcal{O}_B (12)$.
Note also that when 
$V_2 \simeq 
\mathcal{O}_B (5) \oplus 
\mathcal{O}_B (4)^{\oplus 2} \oplus 
\mathcal{O}_B (6)^{\oplus 3}$, we have 
\[
\mathbb{S}^2 (V_2) \simeq 
\mathcal{O}_B (8)^{\oplus 3} \oplus
\mathcal{O}_B (9)^{\oplus 2} \oplus
\mathcal{O}_B (10)^{\oplus 7} \oplus
\mathcal{O}_B (11)^{\oplus 3} \oplus
\mathcal{O}_B (12)^{\oplus 6}.  
\]
Let us denote by 
$\mathcal{O}_B (8)^{\oplus 3} \to \mathbb{S}^2 (\bigwedge^2 V_1)$ 
and 
$\mathbb{S}^2 (V_2) \to \mathcal{O}_B (8)^{\oplus 3}$, 
the natural inclusion and the natural projection, 
respectively, of the direct summand 
$\mathcal{O}_B (8)^{\oplus 3}$,  
of $\mathbb{S}^2 (\bigwedge^2 V_1)$ and 
$\mathbb{S}^2 (V_2)$, respectively. 
Then we see easily (by computation) that 
for general $\sigma_2$, 
the composite 
$\mathcal{O}_B (8)^{\oplus 3} \to \mathcal{O}_B (8)^{\oplus 3}$ 
of the three morphisms 
$\mathcal{O}_B (8)^{\oplus 3} \to \mathbb{S}^2 (\bigwedge^2 V_1)$, 
$\mathbb{S}^2(\sigma_2) \circ \frak{c} : 
\mathbb{S}^2 (\bigwedge^2 V_1) \to \mathbb{S}^2 (V_2)$, 
and 
$\mathbb{S}^2 (V_2) \to \mathcal{O}_B (8)^{\oplus 3}$ 
is a surjection (hence an isomorphism).   
This implies that the image 
$(\mathbb{S}^2(\sigma_2) \circ \frak{c}) 
(\mathcal{O}_B (8)^{\oplus 3})$ is a direct summand 
in $\mathbb{S}^2 (V_2)$. 
Thus we obtain a short exact sequence
\[
0 
\to 
(\mathbb{S}^2 (\bigwedge^2 V_1) )/ 
(\mathcal{O}_B (8)^{\oplus 3}) 
\to 
\mathbb{S}^2 (V_2)/
((\mathbb{S}^2(\sigma_2) \circ \frak{c}) 
(\mathcal{O}_B (8)^{\oplus 3}))
\to
\tilde{V}_4
\to
0,  
\]
where  
$(\mathbb{S}^2 (\bigwedge^2 V_1) )/ 
(\mathcal{O}_B (8)^{\oplus 3}) 
\simeq 
\mathcal{O}_B (10)^{\oplus 2} \oplus \mathcal{O}_B (12)$ and 
\begin{multline}
\mathbb{S}^2 (V_2)/
((\mathbb{S}^2(\sigma_2) \circ \frak{c}) 
(\mathcal{O}_B (8)^{\oplus 3})) 
\simeq  \\
\mathcal{O}_B (9)^{\oplus 2} \oplus
\mathcal{O}_B (10)^{\oplus 7} \oplus
\mathcal{O}_B (11)^{\oplus 3} \oplus
\mathcal{O}_B (12)^{\oplus 6}. \notag
\end{multline}
From this we obtain 
$h^1 (\tilde{V}_4 \otimes (L_4^{\prime})^{\otimes (-1)}) = 0$.  
This implies that 
$h^0 (\tilde{V}_4 \otimes (L_4^{\prime})^{\otimes (-1)}) = 26$ 
and 
$h^1 (\tilde{V}_4 \otimes (L_4^{\prime})^{\otimes (-1)}) = 0$ 
for a general $S^{\flat}$ in our construction. 
Thus under a general fixed extension class, 
the choice of the hypersurfce $X^{\flat}$ depends on 
$25$ parameters. 

Thus we have $15 + 25 = 40$ parameters. 
But we have to subtract 
$\dim \mathrm{Aut} (\mathbb{P}(V_1) / B) = 10$. 
So the surfaces $S^{\flat}$'s in our construction 
fill up a $30$--dimensional stratum in the moduli space 
of regular surfaces with $c_1^2 = 8$ and $p_g = 4$. 
Thus we obtain the following:
\begin{proposition}      \label{prop:m0flat}
In the moduli space of minimal regular surfaces with 
$c_1^2 = 8$ and $p_g = 4$, there exists a $30$--dimensional 
stratum $\mathcal{M}_0^{\flat}$ whose members $S^{\flat}$'s 
admit non--hyperelliptic genus $3$ fibrations 
$f^{\flat} : S^{\flat} \to B$'s with 
$V_1 \simeq 
\mathcal{O}_B (1) \oplus 
\mathcal{O}_B (3) \oplus 
\mathcal{O}_B (3)$ 
and 
$L_4^{\prime} \simeq \mathcal{O}_B (10)$, 
and such that the members $S^{\flat}$'s have canonical maps 
of degree $4$ and no canonical involution.  
\end{proposition}

\begin{remark}         \label{rm:refmflat}
Minimal regular surfaces with $p_g = 4$ 
having canonical maps of degree $4$ can be found, e.g.,   
in \cite{gpclgalloisqd-1}, \cite{gpclgalloisqd-2}, 
and \cite{caninvpg4c8}. 
Note however that our general $S^{\flat}$ has no canonical 
involution, hence different from these examples. 
Note also that our $S^{\flat}$ is topologically simply connected, 
hence different from those in \cite{onbicanonicalmaps}. 
The author does not know any reference on minimal regular 
surfaces with $c_1^2 = 8$, $p_g = 4$ and no canonical 
involution having canonical maps of degree $4$.  
See, e.g., \cite{catlipign} and \cite{caninvpg4c8} 
for references on surfaces with $c_1^2 = 8$ and $p_g = 4$.  
\end{remark}

\medskip



{\sc General members of $\mathcal{M}_0$}

Let $S$ be a member of the stratum $\mathcal{M}_0$ 
given in Bauer--Pignatelli \cite{caninvpg4c8}.
If $S$ is general, 
then $S$ is the minimal desingularization of 
a double cover of $\mathbb{P}^1 \times \mathbb{P}^1$ 
branched along a member of 
the linear system $|8\varDelta_0 + 10 \varGamma|$, 
where $\varDelta_0$ and $\varGamma$ are fibers of 
the first projection and the 2nd projection, respectively,   
and the branch curve has exactly $8$ quadruple points and 
no other essential singularity (see \cite[Theorem 3.5]{caninvpg4c8}). 
Thus $S$ has a hyperelliptic genus $3$ fibration 
$f: S \to B$ that is induced by the second projection 
$\frak{q} : \mathbb{P}^1 \times \mathbb{P}^1 \to B = \mathbb{P}^1$.     

\begin{lemma}   \label{lm:genmembm0}
Let $S$ be a general member of the stratum $\mathcal{M}_0$, 
and $f : S \to B$, the hyperelliptic genus $3$ fibration 
given as above. 
If the corresponding point $[S] \in \mathcal{M}_0$ is general, 
then all the fibers of $f : S \to B$ are $2$--connected. 
Moreover, 
$V_1 = f_* (\omega_{S | B}) \simeq 
\mathcal{O}_B (2) \oplus 
\mathcal{O}_B (2) \oplus 
\mathcal{O}_B (3)$
and 
$L = \mathrm{ker}\, \sigma_2 \simeq \mathcal{O}_B (2)$ 
hold.
\end{lemma}

Proof. 
As given in \cite[Proposition 8]{notesonI}, there exists a 
member of $\mathcal{M}_0$  such that all the fibers 
of the associated hyperelliptic genus $3$ fibrations 
are $2$--connected
(but such that 
$V_1 \simeq \mathcal{O}_B (1) \oplus \mathcal{O}_B (3) 
\oplus \mathcal{O}_B (3)$). 
Since the condition of all the fibers being $2$--connected 
is an open condition in $\mathcal{M}_0$ (see \cite[Lemma 4.1]{notesonI}), 
it follows that for any general member $S$ of $\mathcal{M}_0$ 
all of the fibers of $f : S \to B$ are $2$--connected. 
Moreover, by \cite[Theorem 1]{notesonI}, we have 
$\deg L = 2 (\chi (\mathcal{O}_S) + 2) - 8 -  c_1^2 (S) / 2 = 2$. 
Thus we only need to prove that 
$V_1 \simeq 
\mathcal{O}_B (2) \oplus 
\mathcal{O}_B (2) \oplus 
\mathcal{O}_B (3)$. 
 
Let $\{ w_1, \ldots , w_8 \}$ be the set of the 
quadruple points of the branch divisor of 
$S \to \mathbb{P}^1 \times \mathbb{P}^1$. 
Since 
$V_1 = f_* (\omega_{S | B}) \simeq 
\frak{q}_* (\mathcal{O}_{\mathbb{P}^1 \times \mathbb{P}^1}
           (2 \varDelta_0 + 5 \varGamma - \sum_{i=1}^8 w_i ))$,   
we have for any general choice of $w_1$, $\cdots$, $w_8$ 
a natural short exact sequence  
\begin{equation}    \label{eql:exactcqwi}
0 \to V_1 
   \to \frak{q}_* \mathcal{O}_{\mathbb{P}^1 \times \mathbb{P}^1}
              (2 \varDelta_0 + 5 \varGamma)
  \to \bigoplus_{i=1}^8 \mathbb{C}_{\frak{q} (w_i)} 
  \to 0,
\end{equation}
where $\mathbb{C}_{\frak{q} (w_i)}$ is the skyscraper sheaf 
supported on $\frak{q} (w_i)$ 
and whose stalk is isomorphic to $\mathbb{C}$. 

Note that we have 
$(\frak{q}_* \mathcal{O}_{\mathbb{P}^1 \times \mathbb{P}^1}
           (2 \varDelta_0 + 5 \varGamma)) \otimes \mathcal{O}_B (-3) 
   \simeq \mathcal{O}_B (2)^{\oplus 3}$. 
Form this we see that, if we take $w_1$, $\ldots$, $w_8$ in general 
position, then 
\[
h^0 ((\frak{q}_* \mathcal{O}_{\mathbb{P}^1 \times \mathbb{P}^1}
           (2 \varDelta_0 + 5 \varGamma - \sum_{i=1}^8 w_i))  
    \otimes \mathcal{O}_B (-3))  = 9 - 8 = 1. 
\]
This implies that for
$w_1$, $\ldots$, $w_8$ in general 
position,
the morphism 
\[
H^0 ((\frak{q}_* \mathcal{O}_{\mathbb{P}^1 \times \mathbb{P}^1}
           (2 \varDelta_0 + 5 \varGamma))  
    \otimes \mathcal{O}_B (-3))
   \to
H^0((\bigoplus_{i=1}^8 \mathbb{C}_{\frak{q} (w_i)})
                  \otimes \mathcal{O}_B (-3)))
\]
associated to the exact sequence 
(\ref{eql:exactcqwi}) tensored by $\mathcal{O}_B (-3)$ 
is surjective.
By this together with 
$(\frak{q}_* \mathcal{O}_{\mathbb{P}^1 \times \mathbb{P}^1}
           (2 \varDelta_0 + 5 \varGamma)) \otimes \mathcal{O}_B (-3) 
   \simeq \mathcal{O}_B (2)^{\oplus 3}$, 
we obtain 
$H^1 (V_1 \otimes \mathcal{O}_B (-3)) = 0$. 
Thus if we write $V_1$ as a direct sum of invertible sheaves, 
each direct summand has degree at least $2$. 
Since $\deg V_1 = 7$, this implies 
$V_1 \simeq \mathcal{O}_B (2) \oplus 
            \mathcal{O}_B (2) \oplus 
            \mathcal{O}_B (3)$, 
hence the assertion.  \qed

%
\medskip

Let $x_0$, $x_1$, and $x_2$ be local bases of the 
direct summands $\mathcal{O}_B (2)$,
$\mathcal{O}_B (2)$, and $\mathcal{O}_B (3)$, 
respectively, of the sheaf 
$V_1  =  
\mathcal{O}_B (2) \oplus \mathcal{O}_B (2) \oplus \mathcal{O}_B (3)$. 
Then as in the case of $\mathcal{M}_0^{\sharp}$, 
the morphism 
$\varPhi_{|\mathcal{O}_{\mathbb{P} (V_1)} (1) 
\otimes \mathrm{pr}_{\mathbb{P} (V_1)}^* ( \omega_B )|} : 
\mathbb{P} (V_1) \to \mathbb{P}^3$
is the contraction of 
$\{ x_2 = 0 \} \simeq 
\mathbb{P} (\mathcal{O}_B (2) \oplus \mathcal{O}_B (2)) 
\subset \mathbb{P} (V_1)$ to a line in $\mathbb{P}^3$. 
Let $C \subset \mathbb{P} (V_1)$ be the relative 
conic cut out by the natural inclusion
$L = \mathrm{ker}\, \sigma_2 \to \mathbb{S}^2 (V_1)$. 
Since the canonical map $\varPhi_{|K_S|}$ of $S$ is 
the composite of the natural projection 
$S \to C \subset \mathbb{P} (V_1)$ and the morphism 
$\varPhi_{|\mathcal{O}_{\mathbb{P} (V_1)} (1) 
\otimes \mathrm{pr}_{\mathbb{P} (V_1)}^* ( \omega_B )|} : 
\mathbb{P} (V_1) \to \mathbb{P}^3$, 
we obtain the following: 

\begin{lemma}
Let $S$ be a general member of the stratum $\mathcal{M}_0$ 
in \cite{caninvpg4c8}.  Then the canonical map $\varPhi_{|K_S|}$ of 
$S$ is a morphism of degree $2$ onto its image. 
\end{lemma}

%

\medskip


{\sc Some special members of $\mathcal{M}_0$}

In \cite[Proposition 8]{notesonI}, 
we have constructed a family of members $S$'s of 
$\mathcal{M}_0$ for which the associated 
hyperelliptic genus $3$ fibrations 
$f: S \to B$'s have 
$V_1 \simeq \mathcal{O}_B (1) \oplus \mathcal{O}_B (3) 
\oplus \mathcal{O}_B (3)$ and 
$L = \mathrm{ker}\, \sigma_2 \simeq \mathcal{O}_B (2)$. 
Let us compute the dimension of the substratum 
filled up by $S$'s  
of this family.

%

We have already shown in \cite[Proposition 8]{notesonI} that for 
general $L \to \mathbb{S}^2 (V_1)$ and 
general $(V_2^-)^{\otimes 2} \to A_4$  
the canonical image $\varPhi_{|K_S|} (S)$ of $S$ 
is a non--singular quadric 
$\mathcal{Q} \simeq \mathbb{P}^1 \times B \subset \mathbb{P}^3$.  
The hyperelliptic genus $3$ fibration $f: S \to B$ is 
given by the pullback of one of the two rulings of the quadric $\mathcal{Q}$. 
Thus in order to compute the dimension of the stratum in the 
moduli space, we only need to compute the dimension of the 
family of isomorphism classes (up to automorphisms of the base 
curve $B$) of the hyperelliptic genus $3$ fibrations constructed above.   

%
  
Again we shall count the number of parameters used in the construction.   
For the detail of the construction, see \cite{notesonI}. 

First, the choice of the relative conic 
$C \in |\mathcal{O}_{\mathbb{P}(V_1)} (2) \otimes 
\mathrm{pr}_{\mathbb{P}(V_1)}^*L^{\otimes (-1)}|$ 
cut out by $L = \mathrm{ker}\, \sigma_2 \to \mathbb{S}^2 (V_1)$ 
depends on 
$\dim \mathbb{P} (\mathrm{Hom}\, (L,\, \mathbb{S}^2 (V_1)))$ 
parameters.   
Since $L \simeq \mathcal{O}_B (2)$ and 
$\mathbb{S}^2 (V_1) \simeq 
\mathcal{O}_B (2) \oplus \mathcal{O}_B (4)^{\oplus 2}  
\oplus \mathcal{O}_B (6)^{\oplus 3}$, 
we have 
$h^0 (\mathbb{S}^2 (V_1) \otimes L^{\otimes (-1)}) = 22$. 
Thus the choice of $C$ depends on $21$ parameters. 

Next, under a fixed $C$, the choice of 
the branch divisor of $X \to C$ cut out by 
$(V_2^-)^{\otimes 2} \to 
A_4 = 
\mathrm{Cok}\, (\mathbb{S}^2 (V_1) \otimes L \to \mathbb{S}^4 (V_1))$ 
depends on 
$\dim \mathbb{P} (\mathrm{Hom}\, $ $((V_2^-)^{\otimes 2} ,\, A_4 ))$ 
parameters. 
Let us compute this dimension. 
First note that we have the short exact sequence 
$0 \to \mathbb{S}^2 (V_1) \otimes L \to \mathbb{S}^4 (V_1) \to A_4 \to 0$. 
Note also that we have in our case
$\mathbb{S}^2 (V_1) \otimes L \simeq 
\mathcal{O}_B (4) 
\oplus \mathcal{O}_B (6)^{\oplus 2}  
\oplus \mathcal{O}_B (8)^{\oplus 3}$ 
and 
\[
\mathbb{S}^4 (V_1) \simeq 
\mathcal{O}_B (4) 
\oplus \mathcal{O}_B (6)^{\oplus 2}  
\oplus \mathcal{O}_B (8)^{\oplus 3}
\oplus \mathcal{O}_B (10)^{\oplus 4}  
\oplus \mathcal{O}_B (12)^{\oplus 5}. 
\]
Let us denote by 
$\mathbb{S}^4 (V_1) \to
\mathcal{O}_B (4) 
\oplus \mathcal{O}_B (6)^{\oplus 2}  
\oplus \mathcal{O}_B (8)^{\oplus 3}$
the natural projection onto 
the direct summand 
$\mathcal{O}_B (4) 
\oplus \mathcal{O}_B (6)^{\oplus 2}  
\oplus \mathcal{O}_B (8)^{\oplus 3}$
of $\mathbb{S}^4 (V_1)$. 
Then we see easily (by computation) that 
for general 
$L \simeq \mathcal{O}_B (2) \to \mathbb{S}^2 (V_1)$  
the composite 
$\mathbb{S}^2 (V_1) \otimes L \to 
\mathcal{O}_B (4) 
\oplus \mathcal{O}_B (6)^{\oplus 2}  
\oplus \mathcal{O}_B (8)^{\oplus 3}$ 
of the morphism 
$\mathbb{S}^2 (V_1) \otimes L \to \mathbb{S}^4 (V_1)$ 
and the projection 
$\mathbb{S}^4 (V_1) \to
\mathcal{O}_B (4) 
\oplus \mathcal{O}_B (6)^{\oplus 2}  
\oplus \mathcal{O}_B (8)^{\oplus 3}$ is a surjection 
(hence an isomorphism). 
Form this it follows that 
$A_4 \simeq \mathcal{O}_B (10)^{\oplus 4}  
\oplus \mathcal{O}_B (12)^{\oplus 5}$, 
hence 
$h^0 (A_4 \otimes (V_2^-)^{\otimes (-2)}) = 19$. 
Thus under a fixed $C$, the choice of $X$ depends on 
$18$ parameters. 
 
Thus we have $21 + 18 = 39$ parameters. 
We however need to subtract 
$\dim \mathrm{Aut} (\mathbb{P}(V_1) / B) + 
\dim \mathrm{Aut} (B) = 13$.  
So the surfaces $S$'s in our construction fill up a 
$26$--dimensional stratum. 
Thus we obtain the following:   

\begin{proposition}      \label{prop:m0sp}
In the stratum $\mathcal{M}_0$ of Bauer--Pignatelli \cite{caninvpg4c8}, 
there exists a $26$--dimensional 
substratum $\mathcal{M}_0^{\mathrm{sp}}$ whose members $S$'s 
admit hyperelliptic genus $3$ fibrations 
$f : S \to B$'s with 
$V_1 \simeq 
\mathcal{O}_B (1) \oplus 
\mathcal{O}_B (3) \oplus 
\mathcal{O}_B (3)$ 
and 
$L = \mathrm{ker}\, \sigma_2 \simeq \mathcal{O}_B (2)$, 
and such that the members $S$'s have canonical maps of degree $4$.  
\end{proposition}



{\sc Deformation of general members of $\mathcal{M}_0$}

Now let us deform  
general members $S$'s of the stratum $\mathcal{M}_0$ 
to members of the stratum $\mathcal{M}_0^{\sharp}$.  

\begin{lemma}   \label{lm:v2m0}
Let $f : S \to B$ be a fibration as in Lemma \ref{lm:genmembm0}.
Then $V_2^- \simeq \mathcal{O}_B (5)$. 
Moreover, if $[S] \in \mathcal{M}_0$ is general, 
then 
$V_2^+ \simeq \mathcal{O}_B(5)^{\oplus 4}  
\oplus \mathcal{O}_B (6)$.  
\end{lemma} 

Proof. 
First, note that for a general 
$\nu : \mathcal{O}_B \to \mathcal{O}_B (2)^{\oplus 3}$ 
we have 
$\mathrm{Cok}\, \nu \simeq \mathcal{O}_B(3)^{\oplus 2}$. 
Indeed, if $\nu$ is general, 
the cokernel $\mathrm{Cok}\, \nu$ is a locally free sheaf
of degree $6$ that allows a surjection from 
$\mathcal{O}_B (2)^{\oplus 3}$, 
hence 
$\mathrm{Cok}\, \nu \simeq \mathcal{O}_B (3)^{\oplus 2}$ 
or
$\mathrm{Cok}\, \nu \simeq 
\mathcal{O}_B (2) \oplus \mathcal{O}_B (4)$.
So let us denote by 
$\varLambda_1 \subset 
\mathrm{Hom}(\mathcal{O}_B , \mathcal{O}_B(2)^{\oplus 3}) $ 
the set of 
$\nu$'s such that 
$\mathrm{Cok}\, \nu \simeq \mathcal{O}_B (3)^{\oplus 2}$, 
and by 
$\varLambda_2 \subset 
\mathrm{Hom}(\mathcal{O}_B , \mathcal{O}_B(2)^{\oplus 3}) $ 
the set of 
$\nu$'s such that 
$\mathrm{Cok}\, \nu \simeq 
\mathcal{O}_B (2) \oplus \mathcal{O}_B (4)$.   
Then we have 
\[
\dim \mathbb{P} (\varLambda_1) = 
\dim \mathrm{Hom} (\mathcal{O}_B(2)^{\oplus 3},\, 
\mathcal{O}_B (3)^{\oplus 2}) 
-\dim \mathrm{Aut} (\mathcal{O}_B (3)^{\oplus 2}) = 8,   
\]
and, in the same way, 
$\dim \mathbb{P} (\varLambda_1) = 7$. 
Thus we have 
$\mathrm{Cok}\, \nu \simeq \mathcal{O}_B (3)^{\oplus 2}$
for general $\nu$'s. 

Take a general 
$L = \mathcal{O}_B (2) \to 
\mathbb{S}^2 (V_1) \simeq
\mathcal{O}_B (4)^{\oplus 3} \oplus 
\mathcal{O}_B (5)^{\oplus 2} \oplus 
\mathcal{O}_B (6)$. 
Let 
$\nu^{\prime} : L \simeq \mathcal{O}_B (2) 
\to \mathcal{O}_B (4)^{\oplus 3}$ be the 
composite of this morphism and the natural 
projection to the first direct summand 
$\mathcal{O}_B (4)^{\oplus 3}$.
Then by what we have just remarked, 
we have 
$\mathrm{Cok}\, \nu^{\prime} \simeq \mathcal{O}_B (5)^{\oplus 2}$. 
Thus the morphism 
$H^0 (L \otimes \mathcal{O}_B (-6)) \to 
H^0((\mathcal{O}_B (4)^{\oplus 3}) \otimes \mathcal{O}_B (-6))$ 
induced by the morphism 
$L \to \mathcal{O}_B (4)^{\oplus 3}$ is an isomorphism. 
From this we see that the morphism 
$H^1 (L \otimes \mathcal{O}_B (-6)) \to 
H^1 (V_2^+ \otimes \mathcal{O}_B (-6))$ 
induced by the exact sequence 
\[
0 \to L \otimes \mathcal{O}_B (-6) 
  \to \mathbb{S}^2 (V_1) \otimes \mathcal{O}_B (-6) 
  \to V_2^+ \otimes \mathcal{O}_B (-6) \to 0  
\]
is an isomorphism, hence 
$H^1 (V_2^+ \otimes \mathcal{O}_B (-6)) = 0$. 
This implies that if we write $V_2^+$ as 
a direct summand of invertible sheaves, 
every direct summand has degree at least $5$. 
Since $\deg V_2^+ = 26$, it follows that 
$V_2^+ \simeq \mathcal{O}_B (5)^{\oplus 4} 
\oplus \mathcal{O}_B (6)$. 
The assertion $V_2^- \simeq \mathcal{O}_B (5)$ follows 
from \cite[Proposition 4]{notesonI}. \qed

The following Proposition shows that in the 
moduli space of regular surfaces with $c_1^2 = 8$ and $p_g = 4$, 
the stratum $\mathcal{M}_0$ in \cite{caninvpg4c8} 
lies at the boundary of our 
stratum $\mathcal{M}_0^{\sharp}$.  

\begin{proposition}  \label{prop:defgeneral}
Let $\mathcal{M}_0$ be the $28$--dimensional stratum 
given in the classification list in Bauer--Pignatelli \cite{caninvpg4c8} 
of minimal regular surfaces with $c_1^2 = 8$ and $p_g = 4$ 
and with canonical involutions.
Then for any general member $S$ of $\mathcal{M}_0$, 
there exists a deformation family 
$\pi_{\mathcal{S}} : \mathcal{S} \to T$ of 
$S = \pi_{\mathcal{S}}^{-1} (t_0)$ such that 
for any general $t \in T$ 
the fiber $S_t = \pi_{\mathcal{S}}^{-1} (t)$ is a member 
of the $32$--dimensional stratum $\mathcal{M}_0^{\sharp}$    
constructed in Proposition \ref{prop:m0sharp}. 
\end{proposition}

Proof. 
Let $S$ be a general member of 
$\mathcal{M}_0$, and $f : S \to B$, 
its associated hyperelliptic genus $3$ fibration 
as in Lemmas \ref{lm:genmembm0} and \ref{lm:v2m0}. 
We shall deform this fibration $f : S \to B$ using 
Theorem \ref{thm:exstcthm}. 

Let $L = \mathrm{ker}\, \sigma_2 $ be the kernel of the 
multiplication morphism $\sigma_2 : \mathbb{S}^2 (V_1) \to V_2$, 
and 
$C \in 
|\mathcal{O}_{\mathbb{P}(V_1)} (2) \otimes 
\mathrm{pr}_{\mathbb{P}(V_1)}^* L^{\otimes (-1)}|$, 
the relative conic in $\mathbb{P}(V_1)$ determined by 
the natural inclusion $L \to \mathbb{S}(V_1)$. 
Note that for general $[S] \in \mathcal{M}_0$ the 
relative conic $C$ is non--singular. 
Indeed, since we have 
$\mathbb{S}^2 (V_1) \otimes L^{\otimes (-1)} \simeq 
\mathcal{O}_B(2)^{\oplus 3} \oplus 
\mathcal{O}_B(3)^{\oplus 2} \oplus  
\mathcal{O}_B (4)$ 
by Lemma \ref{lm:genmembm0}, the linear system 
$|\mathcal{O}_{\mathbb{P}(V_1)} (2) \otimes 
\mathrm{pr}_{\mathbb{P}(V_1)}^* L^{\otimes (-1)}|$ 
is free from base points. 
So the smoothness of $C$ follows from Bertini's Theorem 
(note here that taking general $[S] \in \mathcal{M}_0$ corresponds 
to taking general $L \to \mathbb{S}^2 (V_1)$ and general 
$L_4^{\prime} = (V_2^-)^{\otimes 2} \to A_4$). 

Let us check the conditions in Theorem \ref{thm:exstcthm} for our case.   
The condition 1) immediately follows from Lemmas 
\ref{lm:genmembm0} and \ref{lm:v2m0}, because 
by these lemmas we have 
$\mathbb{S}^2 (V_1) \simeq \mathcal{O}_B(4)^{\oplus 3}
\oplus \mathcal{O}_B(5)^{\oplus 2} 
\oplus \mathcal{O}_B (6)$, 
$V_2^+ \simeq \mathcal{O}_B(5)^{\oplus 4} 
\oplus \mathcal{O}_B (6)$, 
and $V_2^- \simeq \mathcal{O}_B(5)$. 
So we only need to check the conditions 2) and 3).

First let us check the condition 2). For this, recall that we have  
$L_4^{\prime}  = (V_2^-)^{\otimes 2} \simeq \mathcal{O}_B (10)$.
Recall also that we have the natural short exact sequence
$0 \to \mathbb{S}^2 (\bigwedge^2 V_1) 
  \to \mathbb{S}^2 (V_2) 
  \to \tilde{V}_4 \to 0$.    
Since we have $V_2 \simeq \mathcal{O}_B(5)^{\oplus 5} 
\oplus \mathcal{O}_B (6)$ 
by Lemmas \ref{lm:genmembm0} and \ref{lm:v2m0}, 
we see from this short excrescence 
tensored by $(L_4^{\prime})^{\otimes (-1)}$ that 
$h^1 (\tilde{V}_4 \otimes (L_4^{\prime})^{\otimes (-1)}) = 0$.  
Thus we have the condition 2) for our general $[S] \in \mathcal{M}_0$. 

Next, let us check the condition 3). 
Recall that we have 
$V_2 = V_2^+ \oplus V_2^-$ and  
$\mathbb{S}^2 (V_2) \simeq \mathbb{S}^2 (V_2^+) \oplus 
(V_2^+ \otimes V_2^-) \oplus (V_2^-)^{\otimes 2}$. 
Meanwhile by Lemma \ref{lm:v2m0}, we have 
\begin{align}
\mathbb{S}^2 (V_2^+) &\simeq 
\mathcal{O}_B(10)^{\oplus 10} \oplus 
\mathcal{O}_B(11)^{\oplus 4} \oplus 
\mathcal{O}_B (12), \notag \\
(V_2^+ \otimes V_2^-) &\simeq 
\mathcal{O}_B(10)^{\oplus 4} \oplus 
\mathcal{O}_B (11),  \notag \\
(V_2^-)^{\otimes 2} &\simeq L_4^{\prime} \simeq 
\mathcal{O}_B (10).  \notag
\end{align}
Thus the linear subsystem of 
$|\mathcal{O}_{\mathbb{P} (V_2)} (2) \otimes 
\mathrm{pr}_{\mathbb{P} (V_2)}^* (L_4^{\prime})^{\otimes (-1)}|$ 
given by the linear subspace generated by   
\begin{multline}
\{ (-\delta ,\, 0, \, \mathrm{id}_{(V_2^-)^{\otimes 2}}) : 
\delta \in 
\mathrm{Hom} ((V_2^-)^{\otimes 2}, \, \mathbb{S}^2 (V_2^+)) \} \\
\subset H^0 (\mathcal{O}_{\mathbb{P} (V_2)} (2) \otimes 
\mathrm{pr}_{\mathbb{P} (V_2)}^* (L_4^{\prime})^{\otimes (-1)}) \notag
\end{multline}
is base point free on $\mathbb{P} (V_2)$, where 
$(-\delta ,\, 0, \, \mathrm{id}_{(V_2^-)^{\otimes 2}})$ 
is regarded as an element of 
$H^0 (\mathcal{O}_{\mathbb{P} (V_2)} (2) \otimes 
\mathrm{pr}_{\mathbb{P} (V_2)}^* (L_4^{\prime})^{\otimes (-1)})$ through the 
natural isomorphism 
\begin{multline}
H^0 (\mathcal{O}_{\mathbb{P} (V_2)} (2) \otimes 
\mathrm{pr}_{\mathbb{P} (V_2)}^* (L_4^{\prime})^{\otimes (-1)})   
\simeq 
\mathrm{Hom} ((V_2^-)^{\otimes 2}, \, \mathbb{S}^2 (V_2))   \\
\simeq 
\mathrm{Hom} ((V_2^-)^{\otimes 2}, \, \mathbb{S}^2 (V_2^+) \oplus 
(V_2^+ \otimes V_2^-) \oplus (V_2^-)^{\otimes 2}).  \notag 
\end{multline} 
Thus, by Bertini's Theorem, for a general 
$\delta \in \mathrm{Hom} ((V_2^-)^{\otimes 2}, \, \mathbb{S}^2 (V_2^+))$, 
the hypersurface $X$ cut out by 
$(-\delta ,\, 0, \, \mathrm{id}_{(V_2^-)^{\otimes 2}})$ in 
$\mathcal{W} |_{B_{t_0}}$ is non-singular, 
where  $\mathcal{W} \subset \mathbb{P} (\mathcal{V}_2)$ 
is the image of a rational veronese map 
as in the proof of Theorem \ref{thm:exstcthm}: 
since we have already checked the conditions  
1) and 2), 
we can take a lifting $\bar{\sigma_2}$ of $\sigma_2$ 
as in the proof of Theorem \ref{thm:exstcthm}, 
and obtain a $\mathcal{W}$ 
associated to this $\bar{\sigma_2} : \mathbb{S}^2 (\mathcal{V}_1) 
\to \mathcal{V}_2$   
(see the proof of Theorem \ref{thm:exstcthm}). 
But this $X$ is the relative canonical model of $S$ 
associated to our general choice of the extension class 
and $\delta$.  Since taking general $[S] \in \mathcal{M}_0$ 
corresponds to taking general extension class 
$0 \to L \to \mathbb{S}^2 (V_1) \to V_2^+ \to 0$ and 
general 
$\delta \in \mathrm{Hom} ((V_2^-)^{\otimes 2},\, A_4)$, 
we have the condition 3) 
for our general $[S] \in \mathcal{M}_0$.  

Thus for any general member $S$ of $\mathcal{M}_0$, 
all the conditions in Theorem \ref{thm:exstcthm} are satisfied. 
It follows that for any general member $S$ of $\mathcal{M}_0$ 
there exists a deformation family 
$(\mathcal{S},\, \mathcal{B},\, T, t_0, 
\frak{f}, \, \pi_{\mathcal{S}}, \pi_{\mathcal{B}})$ 
of $f: S \to B$ as in the proof of Theorem \ref{thm:exstcthm}.  
Note that for our case, the effective divisor $\tau$ 
in the proof of Theorem \ref{thm:exstcthm} can be 
taken in such a way that $|\mathrm{supp}\, \tau| = 3$, 
because the natural morphism 
$\mathrm{Hom} (\mathbb{S}^2 (V_1),\, V_2^- )
\to \mathrm{Hom} (L,\, V_2^-)$ is surjective 
(here $|\mathrm{supp}\, \tau|$ denotes the cardinality 
of the support of the divisor $\tau$). 
Thus, from the construction of the deformation family 
in the proof of Theorem \ref{thm:exstcthm}, 
we see easily that for any general $t \in T$ the 
fiber $S_t = \pi_{\mathcal{S}}^{(-1)} (t)$ 
is a member of the $32$--dimensional family $\mathcal{M}_0^{\sharp}$. 
\qed  


\bigskip

{\sc Deformation of general members of $\mathcal{M}_0^{\mathrm{sp}}$}

Finally Let us deform general members $S$'s of 
the substratum $\mathcal{M}_0^{\mathrm{sp}}$. 

\begin{lemma}  \label{lm:spv2a4isomcls}
Let $S$ be a member of the substratum $\mathcal{M}_0^{\mathrm{sp}}$, 
and $f: S \to B$, the hyperelliptic genus $3$ fibration 
given in Proposition \ref{prop:m0sp}. 
Then $V_2^- \simeq \mathcal{O}_B(5)$. 
Moreover if $[S] \in \mathcal{M}_0^{\mathrm{sp}}$ is general, 
then 
$V_2^+ \simeq 
\mathcal{O}_B(4)^{\oplus 2} \oplus \mathcal{O}_B(6)^{\oplus 3}$ 
and 
$A_4 \simeq 
\mathcal{O}_B(10)^{\oplus 4} \oplus \mathcal{O}_B(12)^{\oplus 5}$. 
\end{lemma}

Proof. 
By Lemma \ref{lm:v2minus}, 
we have 
$V_2^- \simeq (\det V_1) \otimes L^{\otimes (-1)} \simeq \mathcal{O}_B (5)$. 
Now note that 
we have the short exact sequence 
$0 \to L \to \mathbb{S}^2 (V_1) \to V_2^+ \to 0$ 
(see \cite[Lemma 3.4]{notesonI}), 
and that 
$L \simeq \mathcal{O}_B (2)$ and 
$ \mathbb{S}^2 (V_1) \simeq 
\mathcal{O}_B(2) \oplus 
\mathcal{O}_B(4)^{\oplus 2} \oplus 
\mathcal{O}_B(6)^{\oplus 3}$. 
So if $L \to \mathbb{S}^2 (V_1)$ is general 
then the composite 
$L \simeq \mathcal{O}_B(2) \to \mathcal{O}_B(2)$ 
of the morphism $L \to \mathbb{S}^2 (V_1)$ and the 
natural projection $\mathbb{S}^2 (V_1) \to \mathcal{O}_B(2)$ 
is an isomorphism, hence 
$V_2^+ \simeq \mathrm{Cok}\, (L \to \mathbb{S}^2 (V_1)) 
\simeq \mathcal{O}_B(4)^{\oplus 2} \oplus 
\mathcal{O}_B(6)^{\oplus 3}$. 
As for the assertion $A_4 \simeq 
\mathcal{O}_B(10)^{\oplus 4} \oplus \mathcal{O}_B(12)^{\oplus 5}$, 
we have already shown it 
in the proof of Proposition \ref{prop:m0sp}. \qed

The following proposition shows that in the moduli space 
of minimal regular surfaces with $c_1^2 = 8$ and $p_g = 4$, 
the substratum $\mathcal{M}_0^{\mathrm{sp}}$ lies at the 
boundary of the stratum $\mathcal{M}_0^{\flat}$. 

\begin{proposition}  \label{prop:defspecial}
Let $\mathcal{M}_0^{\mathrm{sp}}$ be the $26$--dimensional substratum 
given in Proposition \ref{prop:m0sp}.
Then for any general member $S$ of $\mathcal{M}_0^{\mathrm{sp}}$, 
there exists a deformation family 
$\pi_{\mathcal{S}} : \mathcal{S} \to T$ of 
$S = \pi_{\mathcal{S}}^{-1} (t_0)$ such that 
for any general $t \in T$ 
the fiber $S_t = \pi_{\mathcal{S}}^{-1} (t)$ is a member 
of the $30$--dimensional stratum $\mathcal{M}_0^{\flat}$    
constructed in Proposition \ref{prop:m0flat}. 
\end{proposition}

Proof. 
We use the same method as 
in the proof of Proposition \ref{prop:defgeneral}. 
Let $S$ be a general member of the substratum $\mathcal{M}_0^{\mathrm{sp}}$, 
and $f: S \to B$, the associated hyperelliptic genus $3$ fibration. 
Let us check the conditions in Theorem \ref{thm:exstcthm}. 

By Lemma \ref{lm:spv2a4isomcls}, 
we have 
$\dim \mathrm{Hom}\, (\mathbb{S}^2 (V_1),\, V_2^-) = 8$ 
and 
$\dim \mathrm{Hom}\, (V_2^+,\, V_2^-)$ $= 4$. 
Thus the condition 1) is satisfied.

Also, since 
$\tilde{V}_4 \simeq 
A_4 \oplus (V_2^+ \otimes V_2^-) \oplus (V_2^-)^{\otimes 2}$ 
by Lemma \ref{lm:v4isomoplus}, 
it follows from Lemma \ref{lm:spv2a4isomcls} that  
$h^1 (\tilde{V}_4 \otimes (L_4^{\prime})^{\otimes (-1)}) = 0$, 
where $L_4^{\prime} = (V_2^-)^{\otimes 2}$.  
Thus the condition 2) is satisfied. 

Now note that our case of 
$V_1 \simeq \mathcal{O}_B (1) \oplus 
\mathcal{O}_B (3) \oplus \mathcal{O}_B (3)$ and 
$L \simeq \mathcal{O}_B (2)$ belongs to  
Case B) of \cite[Proposition 7]{notesonI}.
In the proof for Case B) of \cite[Proposition 7]{notesonI}, 
we have shown that if our $S$ is general, 
then both the associated relative conic 
$C \in 
|\mathcal{O}_{\mathbb{P}(V_1)} (2) \otimes 
\mathrm{pr}_{\mathbb{P}(V_1)}^* L^{\otimes (-1)}|$ 
and the branch divisor of 
the natural double cover $X \to C$ are 
non--singular. 
It follows that the relative minimal model $X$ of our $S$ is 
non--singular. 
Thus the condition 3) is also satisfied.  

Finally note that the morphism 
$\mathrm{Hom}\, (\mathbb{S}^2 (V_1),\, V_2^-) 
\to \mathrm{Hom}\, (L ,\, V_2^-)$ 
induced by the exact sequence 
$0 \to L \to \mathbb{S}^2 (V_1) \to V_2^+ \to 0$ 
is surjective. Thus for our case, we can take the 
divisor $\tau$ in the proof of Theorem \ref{thm:exstcthm} 
in such a way that it is supported on three distinct points
in $B$. 
Then by the same argument as in the proof of 
Proposition \ref{prop:defgeneral}, we obtain the assertion. \qed  




\begin{thebibliography}{11}

\bibitem{pg4c7bauer}
{\sc Bauer,~I.} 
Surfaces with $K^2=7$ and $p_g = 4$, 
{\em Mem.\, Amer.\, Math.\, Soc.}, {\bf 152} (2001), no.\, 721.

\bibitem{bcpmodulic6}
{\sc Bauer,~I., Catanese,~F., Pignatelli,~R.} 
The moduli space of surfaces with $K^2=6$ and $p_g=4$, 
{\em Math. Ann.}, {\bf 336} (2006), no.\, 2, 421--438.


\bibitem{caninvpg4c8}
{\sc Bauer,~I., Pignatelli,~R.} 
Surfaces with $K^2=8$, $p_g = 4$ and canonical involution, 
{\em Osaka J. Math.}, {\bf 46} (2009), no.\, 3, 799--820.


\bibitem{catlipign}
{\sc Catanese,~F., Liu,~W., Pignatelli,~R.} 
The moduli space of even surfaces of general type with $K^2 = 8$, $p_g = 4$ and $q = 0$, 
{\tt arXiv:1209.0034 [math.AG]}, {\bf } (2012),.

\bibitem{fibrationsI'}
{\sc Catanese,~F., Pignatelli,~R.} 
Fibrations of low genus, I, 
{\em Ann. Sci. \'{E}cole Norm. Sup. (4)}, {\bf 39} (2006), 1011--1049.

\bibitem{cilcansfpg4}
{\sc Ciliberto,~C.} 
Canonical surfaces with $p_g=p_a=4$ and $K^2=5, \ldots , 10$, 
{\em Duke Math. J.}, {\bf 48} (1981), no.\, 1, 121--157.

\bibitem{onbicanonicalmaps}
{\sc Ciliberto,~C., Francia,~P., Mendes Lopes,~M.} 
Remarks on the bicanonical map for surfaces of general type, 
{\em Math. Z.}, {\bf 224} (1997), 137--166.

\bibitem{enriquessuperfici}
{\sc Enriques,~F.} 
Le superficie algebriche, 
{\em Nicola Zanichelli, Bologna}, {\bf } (1949).


\bibitem{gpclgalloisqd-2}
{\sc Gallego,~F.J., Purnaprajna,~B.P.} 
Classification of quadruple Galois canonical covers. II., 
{\em J. Algebra}, {\bf 312} (2007), no. 2, 798--828.


\bibitem{gpclgalloisqd-1}
{\sc Gallego,~F.J., Purnaprajna,~B.P.} 
Classification of quadruple Galois canonical covers. I, 
{\em Trans. Amer. Math. Soc.}, {\bf 360} (2008), no. 10, 5489--5507.


\bibitem{smallc1-3}
{\sc Horikawa,~E.} 
Algebraic surfaces of general type with small $c_1^2$. III, 
{\em Invent. Math.}, {\bf 47} (1978), no.\ 3, 209--248.

\bibitem{notesonI}
{\sc Murakami,~M.} 
Notes on hyperelliptic fibrations of genus $3$, I, 
{\tt arXiv:1209.6278 [math.AG]}, {\bf } (2012).

\bibitem{oliverioevenpg4c8}
{\sc Oliverio,~P.} 
On even surfaces of general type with $K^2=8$, $p_g=4$, $q=0$, 
{\em Rend.\, Sem.\, Mat.\, Univ.\, Padova}, {\bf 113} (2005), 1--14.

\bibitem{supinopg4c8}
{\sc Supino,~P.} 
On moduli of regular surfaces with $K^2 = 8$ and $p_g = 4$, 
{\em Port.\, Math.\, (N.S.)}, {\bf 60} (2003) no.\, 3, 353--358.



\end{thebibliography}

{\sc 
Masaaki Murakami, 

University of Bayreuth, Lehrstuhl Mathematik VIII, 

Universitaetsstrasse 30, 
D-95447 Bayreuth, Germany}

{\it E-mail address}:\, \texttt{Masaaki.Murakami@uni-bayreuth.de}

\end{document}